\documentclass[11pt]{article}
\usepackage{amssymb, amsmath, latexsym}

\newtheorem{Theo}{Theorem}[section]
\newtheorem{Lem}{Lemma}[section]
\newtheorem{Pro}{Proposition}[section]
\newtheorem{Rem}{Remark}[section]

\newtheorem{Def}{Definition}[section]

\begin{document}

\title{Fractional SPDEs driven by spatially correlated noise: existence of
the solution and smoothness of its density}
\date{ \ \ }
\author{}
\maketitle

\vspace{-2.5cm}

\begin{center}
{\small \sc Lahcen Boulanba} \\
{\small Universit\'{e} Cadi Ayyad, Facult\'{e} des Sciences Semlalia, D\'{e}%
partement de Math\'{e}matiques} \\
{\small B.P. 2390 Marrakech, 40.000, Maroc, l.boulanba@ucam.ac.ma}\\

\medskip

{\small \sc M'hamed Eddahbi\footnote[1]{The second author would like
to express his gratitude for the opportunity to visit Institut de
Math\'ematiques et de Mod\'elisation de Montpellier (I3M),
University Montpellier 2 where this research was done. He would like
to thank also R\'egion Languedoc--Roussillon for his
financial support.}}\\
{\small D\'{e}partement de Math. \& Info. FSTG, Universit\'{e} Cadi Ayyad}\\
{\small B.P. 549, Marrakech, Maroc, eddahbi@fstg-marrakech.ac.ma} \\

\medskip

{\small \sc Mohamed Mellouk\footnote[2]{Corresponding author}}\\
{\small Universit\'{e} Montpellier 2, Institut de Math\'{e}matiques et de Mod\'{e}lisation}\\
{\small de Montpellier, F-34095 Montpellier Cedex 5, France,
mellouk@math.univ-montp2.fr}
\end{center}

\medskip

\begin{abstract}
In this paper we study a class of stochastic partial differential
equations in the whole space $\mathbb{R}^{d}$, with arbitrary
dimension $d\geq 1$, driven by a Gaussian noise white in time and
correlated in space. The differential operator is a fractional
derivative operator. We show the existence, uniqueness and
H\"{o}lder's regularity of the solution. Then by means of Malliavin
calculus, we prove that the law of the solution has a smooth density
with respect to the Lebesgue measure.

\vspace{2cm}

\noindent \textbf{Key words:} Fractional derivative operator,
stochastic partial differential equation, correlated Gaussian noise,
Fourier transform, Malliavin calculus.

\medskip

\noindent \textbf{AMS Subject Classification (2000):} Primary:
60H15; Secondary: 35R60.
\end{abstract}

\newpage

\section{Introduction {and general framework}}

In this paper we consider the following stochastic partial differential
equation (SPDE for abbreviation) given by
\begin{equation}
\left\{
\begin{array}{l}
\dfrac{\partial u}{\partial t}(t,x)=\mathcal{D}_{\delta }^{\alpha
}u(t,x)+b(u(t,x))+\sigma (u(t,x)){\dot{F}}(t,x), \\
u(0,x)=u_{0}(x),%
\end{array}%
\right.  \label{eq0}
\end{equation}%
where $(t,x)\in [0,T]\times \mathbb{R}^{d}$, $d\geq 1$, $\alpha
=(\alpha _{1},\ldots ,\alpha _{d})$, $\delta =(\delta _{1},\ldots
,\delta _{d})$ and $\mathcal{D}_{\delta }^{\alpha }$ denotes a
non--local fractional differential operator to be defined below.
More precisely, we will assume along this paper that
\begin{equation*}
\alpha _{i}\in ]0,2]\setminus \{1\}\text{ and }|\delta _{i}|\leq \min
\{\alpha _{i},2-\alpha _{i}\}\text{, }i=1,\ldots ,d.
\end{equation*}%
The noise $F(t,x)$ is a martingale measure (in the sense given by Walsh in
\cite{Wa}) to be defined with more details in the sequel. The coefficients $%
b $ and $\sigma :\mathbb{R}\rightarrow \mathbb{R}$ are given functions. We
shall refer to Eq. (\ref{eq0}) as $Eq_{\delta }^{\alpha }(d,b,\sigma )$.

Here, we are interested in solutions which are real valued stochastic
processes. Firstly, we establish sufficient conditions ensuring the
existence, the uniqueness and the H\"{o}lder regularity of solutions $u(t,x)$
whenever they exist. Secondly, we study the existence and smoothness of the
density of the law of the solution $u(t,x)$ for fixed $t>0$ and $x\in
\mathbb{R}^{d}$.\newline
Our results of existence, uniqueness and regularity of the solution extend
those obtained in \cite{DeDo} for $Eq_{\delta }^{\alpha }(1,b,\sigma )$ and
\cite{AM} for $Eq_{0}^{\alpha }(1,b,\sigma )$ both driven by a space--time
white noise $W$ and $1<\alpha <2$. However in their framework, the
regularity of the noise rises up $\alpha >1$ as a necessary condition for
the existence of the stochastic integral as an $L^{2}(\Omega )$ random
variable. This problem does not appear in our case ($\alpha \in
]0,2]\setminus \{1\}$) because the noise $F$ is smoother than $W$. \newline
The equation $Eq_{\delta }^{\alpha }(d,b,\sigma )$ recovers for instance the
stochastic heat equation in spatial dimension $d\geq 1$ studied by many
authors among others \cite{Da}, \cite{Sa}. With the notations adopted above
it corresponds to $Eq_{0}^{2}(d,b,\sigma )$. Indeed, when $\delta _{i}=0$
and $\alpha _{i}=2$ for $i=1,\ldots ,d$ the operator $\mathcal{D}_{0}^{2}$
coincides with the classical Laplacian operator in $\mathbb{R}^{d}$.\newline
Various physical phenomena involving diffusion and interaction of particles
can be described by the equations $Eq_{\delta }^{\alpha }(d,b,\sigma )$ when
suitable assumptions are made on the coefficient $b$ and $\sigma $. The
non--local property in this equation is due to the presence of $\mathcal{D}%
_{\delta }^{\alpha }$ and the non--linearity comes from the general form of $%
b$ and $\sigma $. This equations can also be interpreted as random
perturbation of deterministic equations $Eq_{\delta }^{\alpha }(d,b,0)$ by
multiplicative noise of the form $\sigma (u(t,x)){\dot{F}}(t,x)$.\newline
In probabilistic terms, replacing the Laplacian by its fractional power
(which is an integro--differential operator) leads to interesting and
largely open questions of extensions of results for Brownian motion driven
stochastic equations to those driven by L\'{e}vy stable processes. In the
physical literature, such fractal anomalous diffusions have been recently
enthusiastically embraced by a slew of investigators in the context of
hydrodynamics, acoustics, trapping effects in surface diffusion, statistical
mechanics, relaxation phenomena, and biology (see e.g. \cite{SZF}, \cite{Su}%
, \cite{Z}, \cite{ZA}, \cite{SW}, \cite{W1}).\newline
A probabilistic approach to the equation $Eq_{0}^{\alpha }(d,b,0)$ is made
by means of the Feynman--Kac formula (see \cite{BLW0}, \cite{BLW}, \cite{JMW}%
). Solutions to other particular fractional differential equations are given
as functionals of stable subordinators. This representation provides
explicit form to the density of the $3/2$--stable law and to the density of
escaping island vicinity in vortex medium, see \cite{De}. In other words the
laws of stable L\'{e}vy processes or stable subordinators satisfies
fractional equations like $Eq_{\delta }^{\alpha }(d,b,0)$.

SPDE lies at the intersection of several different areas, and different
groups have taken the subject in their own directions. For example, one
could start with the functional analytic approach to PDE, and develop a
similar framework for SPDE. Reasoning this way, one views solutions of SPDE
as random variables taking values in Sobolev spaces (or other spaces), and
the terms in the equations as operators and operator--valued random
variables.\newline
Another possibility is to think of an SPDE as a stochastic process. In this
approach, we focus on the intriguing properties of Gaussian processes and
martingale measures. We will call this point of view the Walsh approach.%
\newline
To give a flavor of this approach, let us consider two particular processes,
the first one is the superprocess or Dawson--Watanabe process. Roughly
speaking, the density of the its law of this solves the SPDE
\begin{equation*}
\frac{\partial u}{\partial t}(t,x)=\mathcal{D}_{0}^{2}u(t,x)+\sqrt{u(t,x)}{%
\dot{F}}(t,x)\ \ \ \ \ (d=1,\ \mathcal{D}_{0}^{2}=\Delta )
\end{equation*}%
and the superprocess considered in \cite{LMX} whose associated SPDE is given
by
\begin{equation*}
\frac{\partial u}{\partial t}(t,x)=Lu(t,x)-\frac{\partial (gu)}{\partial x}%
(t,x)\dot{W}_{t}+\sqrt{u(t,x)}{\dot{F}}(t,x)
\end{equation*}%
where $L$ is a second order differential operator, $g$ is a given function
and $W$ is a Brownian motion.\newline
Then super--processes arise as a limiting process in population dynamics and
genetics. This connection leads to an analysis of certain SPDE in terms of
particles. One of the sources of SPDE was the Zakai equation of filtering
theory, and there was some work on the stochastic Navier--Stokes equation,
but the theory of SPDE itself was not closely tied to applications. Luckily,
physicists and engineers have now caught up, and SPDE is appearing in more
and more scientific fields. We may mention some of these connections, such
as randomly growing surfaces in solid--state physics and the theory of
randomly moving polymers. Also, the stochastic Navier--Stokes equation has
received much more intensive study in recent years.

Our own motivation was, however, the connection between fractional operators
(fractional power of the Laplacian) and super--processes. For example
equation of the form $Eq_{0}^{\alpha }(d,b,\sigma )$ are related to some
measure--valued processes. More precisely $Eq_{0}^{\alpha }(1,0,\sqrt{\cdot }%
)$ is satisfied by the density for the law of a measure--valued
branching diffusion for $1<\alpha \leq 2$ and $\delta =0$. Another
example comes from Fleming--Viot models, that is, under some
conditions on this diffusion process its law is absolutely
continuous w.r.t. the Lebesgue measure with continuous density
$u(t,x)$ for each $(t,x)\in ]0,+\infty [\times \mathbb{R}$. Moreover
it satisfies the following SPDE
\begin{equation*}
\dfrac{\partial u}{\partial t}(t,x)=\mathcal{D}_{0}^{\alpha }u(t,x)+\sqrt{%
u(t,x)}{\dot{F}}(t,x)-u(t,x)\int_{\mathbb{R}}\sqrt{u(t,y)}{\dot{F}}(t,y)dy,
\end{equation*}%
for $d=1$, $1<\alpha \leq 2$ and ${\dot{F}}(t,x)$ is a space--time white
noise in $\mathbb{R}_{+}\times \mathbb{R}$. For complete description of
measure--valued branching and Fleming--Viot diffusions the reader may
consult \cite{KoSh} (see also \cite{FlVi} and \cite{RoCo}). The above
equation is the perturbation of $Eq_{0}^{\alpha }(1,0,\sqrt{\cdot })$ by $%
u(t,x)\int_{\mathbb{R}}(u(t,y))^{1/2}{\dot{F}}(t,y)dy$. To recover this
equation in $Eq_{\delta }^{\alpha }(d,b,\sigma )$ we allow $\sigma $ to
depend on the whole trajectory of $u$. Nevertheless we are not studying this
problem in this paper and restrict our selves to Lipschitz coefficients. The
case where $0<\alpha <1$ was motivated by the work \cite{Chen} where he
attempts to use the $2/3$--order fractional Laplacian modeling of enhanced
diffusing movements of random turbulent particle resulting from non--linear
inertial interactions. A combined effect of the inertial interactions and
the molecule Brownian diffusivities is found to be the bi--fractal mechanism
behind multifractal scaling in the inertial range of scales of moderate
Reynolds number turbulence. Accordingly, a stochastic equation is proposed
to describe turbulence intermittency. The $2/3$--order fractional Laplacian
representation is also used to construct a fractional Reynolds equation for
non--linear interactions of fluctuating velocity components, underlying
turbulence space--time fractal structures of Levy $2/3$--stable
distribution. The new perspective of his study is that the fractional
calculus is an effective approach modeling of chaotic fractal phenomena
induced by non--linear interactions.

Fractional derivatives and integrals, usually known as \textit{fractional
calculus}, have many uses and they themselves have arisen from certain
requirements in applications, such as fractional integro--differentiation
which has now become a significant topic in mathematical analysis. It has
applications in various fields namely quantitative biology,
electrochemistry, transport theory, probability and potential theory to
mention a few. We refer the reader for a complete survey on the fractional
integrals and derivatives to \cite{Sam} and \cite{Po} (and the references
therein). In recent years, fractional derivative has long been found to be a
very effective means to describe the anomalous attenuation behaviors, some
parametric seismic wave propagation have been model via the fractional
derivative model. Here will not discus this topic and limit our selves to
the above connections between fractional calculus and mathematical and
physical phenomena.\newline
In the literature, various fractional differential operators are defined
(see \cite{DGV}, \cite{MW}, \cite{Po}). The results of this paper apply to
several of them, such as fractional Laplacian, Nishimoto operator and the
non--self adjoint fractional operator introduced in \cite{Ko} and used in
\cite{KoSh} to the study of stochastic partial differential equation.

Basic notations, definitions and preliminary results of the operator $%
\mathcal{D}_{\delta }^{\alpha }$ and the noise $F$ will be presented in the
following subsequent two subsections. The section 2 is devoted to the
existence and uniqueness and the section 3 deals with the H\"{o}lder
regularity result. The section 4 contains the existence and smoothness of
the density of the law of the solution. In appendix we prove some technical
results which will be used in the proofs.

The value of the constants along this article may change from line to line
and some of the standing parameters are not always indicated.

\subsection{The operator $\mathcal{D}_{\protect\delta }^{\protect\alpha }$}

For the sake of notational simplicity we denote by $D_{\delta }^{\alpha }$
the fractional differential operator in dimension $d=1$ and $\mathcal{D}%
_{\delta }^{\alpha }$ in dimension $d\geq 2$. \newline
Let us give the definition of the operator $D_{\delta }^{\alpha }$ (i.e. in
one space dimension).

\begin{Def}
The fractional differential operator $D_{\delta }^{\alpha }$ is a non--local
operator defined via its Fourier transform $\mathcal{F}$ by
\begin{equation*}
\mathcal{F}\left( D_{\delta }^{\alpha }\varphi \right) (\xi )=-|\xi
|^{\alpha }\exp \left( -\imath \delta \frac{\pi }{2}\mathrm{sgn}(\xi
)\right) \mathcal{F}\left( \varphi \right) (\xi ),
\end{equation*}%
where $\imath ^{2}+1=0$.
\end{Def}

The operator $D_{\delta }^{\alpha }$ is a closed, densely defined operator
on $L^{2}(\mathbb{R})$ and it is the infinitesimal generator of a semigroup
which is in general not symmetric and not a contraction. It is self adjoint
only when $\delta =0$ and in this case, it coincides with the fractional
power of the Laplacian.

According to \cite{Ko}, $D_{\delta }^{\alpha }$ can be represented for $%
1<\alpha <2$, by
\begin{equation*}
D_{\delta }^{\alpha }\varphi (x)=\int_{-\infty }^{+\infty }\frac{\varphi
(x+y)-\varphi (x)-y\varphi ^{\prime }(x)}{|y|^{1+\alpha }}\left( \kappa
_{-}^{\delta }\mathbf{1}_{(-\infty ,0)}(y)+\kappa _{+}^{\delta }\mathbf{1}%
_{(0,+\infty )}(y)\right) dy
\end{equation*}%
and for $0<\alpha <1$, by%
\begin{equation*}
D_{\delta }^{\alpha }\varphi (x)=\int_{-\infty }^{+\infty }\frac{\varphi
(x+y)-\varphi (x)}{|y|^{1+\alpha }}\left( \kappa _{-}^{\delta }\mathbf{1}%
_{(-\infty ,0)}(y)+\kappa _{+}^{\delta }\mathbf{1}_{(0,+\infty )}(y)\right)
dy
\end{equation*}%
where $\kappa _{-}^{\delta }$ and $\kappa _{+}^{\delta }$ are two
non--negative constants satisfying $\kappa _{-}^{\delta }+\kappa
_{+}^{\delta }>0$ and $\varphi $ is a smooth function for which the integral
exists, and $\varphi ^{\prime }$ is its derivative. This representation
identifies it as the infinitesimal generator for a non--symmetric $\alpha $%
--stable L\'{e}vy process.

Let $G_{\alpha ,\delta }(t,x)$ denotes the fundamental solution of the
equation $Eq_{\delta }^{\alpha }(1,0,0)$ that is, the unique solution of the
Cauchy problem
\begin{equation*}
\left\{
\begin{array}{l}
\dfrac{\partial u}{\partial t}(t,x)=D_{\delta }^{\alpha }u(t,x), \\
u(0,x)=\delta _{0}(x),\ \ t>0,\text{ }x\in \mathbb{R},%
\end{array}%
\right.
\end{equation*}%
where $\delta _{0}$ is the Dirac distribution. Using Fourier's calculus we
obtain
\begin{equation*}
G_{\alpha ,\delta }(t,x)=\frac{1}{2\pi }\int_{-\infty }^{+\infty }\exp
\left( -\imath zx-t|z|^{\alpha }\exp \left( -\imath \delta \frac{\pi }{2}%
\mathrm{sgn}(z)\right) \right) dz.
\end{equation*}%
The relevant parameters, $\alpha $, called the index of \textit{stability}
and $\delta $ (related to the asymmetry) improperly referred to as the
\textit{skewness} are real numbers satisfying $\alpha \in ]0,2]$ and $%
|\delta |\leq \{\alpha ,2-\alpha \}.$

The function $G_{\alpha ,\delta }(t,\cdot )$ has the following properties.

\begin{Lem}
\label{lgr} For $\alpha \in (0,2]\setminus \{1\}$ such that $|\delta |\leq
\min (\alpha ,2-\alpha )$

\begin{description}
\item[$(i)$] The function $G_{\alpha ,\delta }(t,\cdot )$ is the density of
a L\'{e}vy $\alpha $--stable process in time $t$.

\item[$(ii)$] The function $G_{\alpha ,\delta }(t,x)$ is not in general
symmetric relatively to $x$.

\item[$(iii)$] Semigroup property: $G_{\alpha ,\delta }(t,x)$ satisfies the
Chapman Kolmogorov equation, i.e. for $0<s<t$%
\begin{equation*}
G_{\alpha ,\delta }(t+s,x)=\int_{-\infty }^{+\infty }G_{\alpha ,\delta
}(t,y)G_{\alpha ,\delta }(s,x-y)dy.
\end{equation*}

\item[$(iv)$] Scaling property: $G_{\alpha ,\delta }(t,x)=t^{-1/\alpha
}G_{\alpha ,\delta }(1,t^{-1/\alpha }x)$.

\item[$(v)$] There exists a constant $c_{\alpha }$ such that $0\leq
G_{\alpha ,\delta }(1,x)\leq c_{\alpha }(1+|x|^{1+\alpha })$, for all $x\in
\mathbb{R}$.
\end{description}
\end{Lem}

For $d\geq 1$ and any multi index $\alpha =(\alpha _{1},\ldots ,\alpha _{d})$
and $\delta =(\delta _{1},\ldots ,\delta _{d})$, define the operator $%
\mathcal{D}_{\delta }^{\alpha }$ by
\begin{equation*}
\mathcal{D}_{\delta }^{\alpha }=\sum_{i=1}^{d}D_{\delta _{i}}^{\alpha _{i}},
\end{equation*}%
where $D_{\delta _{i}}^{\alpha _{i}}$ denotes the fractional differential
derivative w.r.t the $i$--th coordinate.

Let $\mathbf{G}_{\alpha ,\delta }(t,x)$ be the Green function of the
deterministic equation $Eq_{\delta }^{\alpha }(d,0,0)$. Clearly
\begin{eqnarray*}
\mathbf{G}_{\alpha ,\delta }(t,x) &=&\prod_{i=1}^{d}G_{\alpha _{i},\delta
_{i}}(t,x_{i}) \\
&=&\frac{1}{(2\pi )^{d}}\int_{\mathbb{R}^{d}}\exp \left( -\imath
\left\langle \xi ,x\right\rangle -t\sum_{i=1}^{d}|\xi _{i}|^{\alpha
_{i}}\exp \left( -\imath \delta _{i}\frac{\pi }{2}\mathrm{sgn}(\xi
_{i})\right) \right) d\xi ,
\end{eqnarray*}%
where $\left\langle \cdot ,\cdot \right\rangle $ stands for the inner
product in $\mathbb{R}^{d}$.

\subsection{The driving noise $F$}

Let $\mathcal{S}(\mathbb{R}^{d+1})$ be the space of Schwartz test functions.
On a complete probability space $(\Omega ,\mathcal{G},P)$, the noise $%
F=\{F(\varphi ),\varphi \in \mathcal{S}(\mathbb{R}^{d+1})\}$ is assumed to
be an $L^{2}(\Omega ,\mathcal{G},P)$--valued Gaussian process with mean zero
and covariance functional given by
\begin{equation*}
J(\varphi ,\psi )=\int_{\mathbb{R}_{+}}ds\int_{\mathbb{R}^{d}}\Gamma
(dx)\left( \varphi (s,\cdot )\ast \widetilde{\psi }(s,\cdot )\right) (x),\ \
\varphi ,\psi \in \mathcal{S}(\mathbb{R}^{d+1}),
\end{equation*}%
where $\widetilde{\psi }(s,x)=\psi (s,-x)$ and $\Gamma $ is a non--negative
and non--negative definite tempered measure, therefore symmetric. Let $\mu $
denote the spectral measure of $\Gamma $, which is also a tempered measure
(see \cite{Sch}) that is $\mu =\mathcal{F}^{-1}(\Gamma )$ and this gives
\begin{equation}
J(\varphi ,\psi )=\int_{\mathbb{R}_{+}}ds\int_{\mathbb{R}^{d}}\mu (d\xi )%
\mathcal{F}\varphi (s,\cdot )(\xi )\overline{\mathcal{F}\psi (s,\cdot )}(\xi
),  \label{cov}
\end{equation}%
where $\overline{z}$ is the complex conjugate of $z.$

\noindent Following the same approach in \cite{Da}, the Gaussian process $F$
can be extended to a worthy martingale measure $M=\{M(t,A):=F([0,t]\times
A)\,:\,t\geq 0,\,A\in \mathcal{B}_{b}(\mathbb{R}^{d})\}$ which shall acts as
integrator, in the sense of Walsh \cite{Wa}, where $\mathcal{B}_{b}(\mathbb{R%
}^{d})$ denotes the bounded Borel subsets of $\mathbb{R}^{d}$. Let $\mathcal{%
G}_{t}$ be the completion of the $\sigma $--field generated by the random
variables $\{M(s,A),\;0\leq s\leq t,\;A\in \mathcal{B}_{b}(\mathbb{R}^{d})\}$%
. The properties of $F$ ensure that the process $M=\{M(t,A),\;t\geq 0,\;A\in
\mathcal{B}_{b}(\mathbb{R}^{d})\}$, is a martingale with respect to the
filtration $\{\mathcal{G}_{t}:t\geq 0\}$.

Then one can give a rigorous meaning to solution of equation $Eq_{\delta
}^{\alpha }(d,b,\sigma )$, by means of a jointly measurable and $\mathcal{G}%
_{t}$--adapted process $\{u(t,x):(t,x)\in \mathbb{R}_{+}\times \mathbb{R}%
^{d}\}$ satisfying, for each $t\geq 0$, a.s. for almost $x\in \mathbb{R}^{d}$%
, the following evolution equation:
\begin{eqnarray}
u(t,x) &=&\int_{\mathbb{R}^{d}}\mathbf{G}_{\alpha ,\delta }(t,x-y)u_{0}(y)dy
\notag \\
&&+\int_{0}^{t}ds\int_{\mathbb{R}^{d}}\mathbf{G}_{\alpha ,\delta
}(t-s,x-y)b(u(s,y))dy  \notag \\
&&+\int_{0}^{t}\int_{\mathbb{R}^{d}}\mathbf{G}_{\alpha ,\delta
}(t-s,x-y)\sigma (u(s,y))M(ds,dy)  \label{eqfrac}
\end{eqnarray}%
Throughout this paper we adopt the following definition.

\begin{Def}
\label{d1}A stochastic process $u$ defined on $\Omega \times \mathbb{R}%
_{+}\times \mathbb{R}^{d}$, which is jointly measurable and $\mathcal{G}_{t}$%
--adapted, is said to be a solution to the fractional SPDE $Eq_{\delta
}^{\alpha }(d,b,\sigma )$, if it is an $\mathbb{R}$--valued fields which
satisfies $(\ref{eqfrac})$ and
\begin{equation*}
\sup_{t\in [0,T]}\sup_{x\in \mathbb{R}^{d}}E|u(t,x)|^{p}<+\infty
\text{, for some }p\geq 2.
\end{equation*}
\end{Def}

Let us formulate our assumption concerning the fractional differential
operator $\mathcal{D}_{\delta }^{\alpha }$ and the correlation of the noise $%
M$.

If we take $\sigma \equiv 1${\ and use the formula (\ref{cov}) the
stochastic integral }%
\begin{equation*}
\int_{0}^{T}\int_{\mathbb{R}^{d}}{\mathbf{G}_{\alpha ,\delta }}%
(T-s,x-y)M(ds,dy)
\end{equation*}%
{appearing in (\ref{eqfrac}) is well defined if and only if }%
\begin{equation*}
\int_{0}^{T}ds\int_{\mathbb{R}^{d}}\mu (d\xi )|{\mathcal{F}\mathbf{G}%
_{\alpha ,\delta }}(s,\cdot )(\xi )|^{2}<+\infty .
\end{equation*}%
{Indeed
\begin{eqnarray*}
&&E\left\vert \int_{0}^{t}\int_{\mathbb{R}^{d}}{\mathbf{G}_{\alpha ,\delta }}%
(t-s,x-y)M(ds,dy)\right\vert ^{2} \\
&=&\int_{0}^{t}ds\int_{\mathbb{R}^{d}}{\Gamma }(dz)({\mathbf{G}_{\alpha
,\delta }}(t-s,\cdot )\ast \widetilde{{\mathbf{G}_{\alpha ,\delta }}}%
(t-s,\cdot ))(x-z) \\
&=&\int_{0}^{t}ds\int_{\mathbb{R}^{d}}\mu (d\xi )|{\mathcal{F}\mathbf{G}%
_{\alpha ,\delta }}(t-s,\cdot )(\xi )|^{2},\ \ ({\mathcal{F}}{\mu }={\Gamma }%
).
\end{eqnarray*}%
For a given multi index }$\alpha =(\alpha _{1},\ldots ,\alpha _{d})$ such
that $0<\alpha _{i}\leq 2$, $i=1,\ldots ,d$ and any $\xi \in \mathbb{R}^{d}$
we use the notation
\begin{equation*}
S_{\alpha }(\xi )=\sum_{i=1}^{d}|\xi _{i}|^{\alpha _{i}}.
\end{equation*}%
The following lemma gives sufficient condition for the existence of
stochastic integrals w.r.t. the martingale measure $M$.

\begin{Lem}
\label{lemmu}There exist two positive constants $c_{1}$ and $c_{2}$ such
that
\begin{equation}
c_{1}\int_{\mathbb{R}^{d}}\frac{\mu (d\xi )}{1+S_{\alpha }(\xi )}\leq
\int_{0}^{T}\int_{\mathbb{R}^{d}}\mu (d\xi )\left\vert \mathcal{F}\mathbf{G}%
_{\alpha ,\delta }(s,\cdot )(\xi )\right\vert ^{2}ds\leq c_{2}\int_{\mathbb{R%
}^{d}}\frac{\mu (d\xi )}{1+S_{\alpha }(\xi )}.  \label{eqfourier}
\end{equation}
\end{Lem}

\noindent \textbf{Proof.} By the definition of $\mathbf{G}_{\alpha ,\delta }$
we have
\begin{equation*}
\mathcal{F}\mathbf{G}_{\alpha ,\delta }(t,\cdot )(\xi )=\exp \left(
-t\sum_{i=1}^{d}|\xi _{i}|^{\alpha _{i}}\exp \left( -\imath \delta _{i}\frac{%
\pi }{2}\mathrm{sgn}(\xi _{i})\right) \right) .
\end{equation*}%
Therefore,
\begin{equation}
\int_{0}^{t}\left\vert \mathcal{F}\mathbf{G}_{\alpha ,\delta }(s,\cdot )(\xi
)\right\vert ^{2}ds=\Phi \left( 2t\sum_{i=1}^{d}|\xi _{i}|^{\alpha _{i}}\cos
\left( \delta _{i}\frac{\pi }{2}\right) \right) ,  \label{fi1}
\end{equation}%
where $\Phi (x)=t(1-e^{-x})/x$. Since $\Phi $ is a decreasing function, we
get
\begin{equation}
\Phi \left( 2tS_{\alpha }(\xi )\right) \leq \Phi \left( 2t\sum_{i=1}^{d}|\xi
_{i}|^{\alpha _{i}}\cos \left( \delta _{i}\frac{\pi }{2}\right) \right) \leq
\Phi \left( 2t\kappa S_{\alpha }(\xi )\right)  \label{fi2}
\end{equation}%
with $\kappa =\min_{1\leq i\leq d}(\cos (\delta _{i}\pi /2))>0$ since $-\pi
/2<\delta _{i}\pi /2<\pi /2$. The mean value theorem applied to the function
$(1+x)(1-e^{-x})$ yields
\begin{equation*}
\frac{t}{1+x}\leq \Phi (x)\leq \frac{2t}{1+x}
\end{equation*}%
for all $x>0$ and $t\in [0,T]$. This together with (\ref{fi1}) and (%
\ref{fi2}) imply that for all $t\in [0,T]$
\begin{equation*}
\frac{t}{1+2tS_{\alpha }(\xi )}\leq \int_{0}^{t}\left\vert \mathcal{F}%
\mathbf{G}_{\alpha ,\delta }(s,\cdot )(\xi )\right\vert ^{2}ds\leq \frac{2t}{%
1+2t\kappa S_{\alpha }(\xi )}.
\end{equation*}%
The conclusion follows by taking the integral over $\mathbb{R}^{d}$ with
respect to the measure ${\mu }$.\hfill $\square $

\begin{Rem}
The upper and lower bounds in $(\ref{eqfourier})$ do not depend on the
parameter $\delta $. When $\alpha _{i}=2$ for all $i=1,\ldots ,d$ then $%
S_{2}(\xi )=\left\Vert \xi \right\Vert _{2}^{2}$ and the bounds in $(\ref%
{eqfourier})$ are the same ones which appeared in \cite{Da} (see also \cite%
{Sa}), that is
\begin{equation*}
\int_{\mathbb{R}^{d}}\frac{\mu (d\xi )}{1+\left\Vert \xi \right\Vert _{2}^{2}%
}<+\infty .
\end{equation*}
\end{Rem}

Our main assumption on the noise relies an the integrability condition
w.r.t. the spectral measure $\mu .$\newline
\noindent \textbf{Assumption }$(\mathbf{H}_{\eta }^{\alpha })$
\begin{equation*}
\int_{\mathbb{R}^{d}}\frac{{\mu }(d\xi )}{\left( 1+S_{\alpha }(\xi )\right)
^{\eta }}<+\infty \ \ \text{where}\ \eta \in ]0,1].
\end{equation*}%
This condition states that there are not too many high frequencies in the
noise, which turns out to be reformulated into a condition on the
integrability of $\Gamma $ (measure characterizing the covariance of the
noise), near the origin in the case where $\Gamma $ is a non--negative
measure on $\mathbb{R}^{d}$. The condition $(\mathbf{H}_{0}^{\alpha })$
means that $\mu $ is a finite measure, which is equivalent to say that $%
\Gamma $ is a uniformly continuous and bounded function on $\mathbb{R}^{d}$.
For this reason we do not consider this case in this paper.\newline
When $\mu $ is the Lebesgue measure on $\mathbb{R}^{d}$ which is the
spectral measure of the white noise on $\mathbb{R}^{d}$ that is the noise
corresponding to $\Gamma =\delta _{0}$, the condition $(\mathbf{H}_{\eta
}^{\alpha })$ is equivalent to $\int_{0}^{+\infty }{r}^{-1+\sum_{i=1}^{d}%
\frac{1}{\alpha _{i}}}\left( 1+r\right) ^{-\eta }dr<+\infty $ which is
finite in and only if $\eta >\sum_{i=1}^{d}(1/\alpha _{i})$.

\subsubsection{Some examples}

Let us now, consider a class of covariances measures $\Gamma $ for which the
condition $(\mathbf{H}_{\eta }^{\alpha })$ gives an optimal criterion. In
this part we assume $\alpha _{i}=2$, $i=1,\ldots ,d$ and use $|x|$ as the
Euclidian norm of $x$.

\begin{description}
\item[1.] Riesz kernels. Let $\Gamma (dx)=|x|^{-\gamma }dx$, $\gamma \in
]0,d[$. It is known that the spectral measure is $\mu (dx)=c|x|^{\gamma
-d}dx $. Then $(\mathbf{H}_{\eta }^{2})$ is satisfied if and only if $\gamma
\in ]0,2\eta \wedge d[$.

\item[2.] Assume that the spectral measure $\mu $ is either finite or
absolutely continuous with respect to Lebesgue measure and $\mu (dx)/dx=\rho
(x)\in L^{p}(\mathbb{R}^{d})$ for some $p\geq 1$. For $\mu $ finite $(%
\mathbf{H}_{\eta }^{2})$ holds for any $\gamma \in ]0,1[$. In the second
case, $(\mathbf{H}_{\eta }^{2})$ holds if $d(p-1)/2p<\eta <1$. That is if $%
p\in [1,d/(d-2)\vee 1]$.

\item[3.] Brownian free field. Let $J(\varphi ,\psi )=\left\langle (-\Delta
+m^{2})^{-1}\varphi ,\psi \right\rangle $, where $\Delta $ is the Laplace
operator in $\mathbb{R}^{d}$ and $m>0$. In this case $\mu (dx)=(2\pi
)^{-d/2}(|x|^{2}+m^{2})^{-1}dx$. Therefore, $(\mathbf{H}_{\eta }^{2})$ holds
if and only if $\eta \in ](d-2)^{+}/2\wedge 1,1[$. That is, for $d=1$, $2$, $%
(\mathbf{H}_{\eta }^{2})$ holds for any $\eta \in (0,1)$, for $d=3$, $(%
\mathbf{H}_{\eta }^{2})$ holds for any $\eta \in ]1/2,1[$; for $d\geq 4$, $(%
\mathbf{H}_{\eta }^{2})$ does not hold.

\item[4.] Bessel Kernel. Assume $\beta >0$ and let $\Gamma $ be the Bessel
kernel of order $\beta $. That is $\Gamma (dx)=B_{\beta }(x)dx$ where $%
B_{\beta }(x)=\nu _{\beta ,d}\int_{0}^{+\infty }t^{-1+(\beta
-d)/2}e^{-t-|x|^{2}/4t}dt$, $x\in \mathbb{R}^{d}$ and $\nu _{\beta ,d}$ is a
constant. The spectral measure $\mu $ of $\Gamma $ is given by $\mathcal{F}%
(B_{\beta })(\xi )=(1+|\xi |^{2})^{-\beta /2}$. Therefore $(\mathbf{H}_{\eta
}^{2})$ is satisfied if and only if
\begin{equation*}
\int_{\mathbb{R}^{d}}(1+|\xi |^{2})^{-\eta -\beta /2}d\xi
=c_{d}\int_{0}^{+\infty }(1+t)^{-\eta -\beta /2}t^{-1+d/2}dt
\end{equation*}%
which is convergent for $d<2\eta +\beta $. Hence $(\mathbf{H}_{\eta }^{2})$
holds for $\eta \in ](d-\beta )^{+}/2,1[$
\end{description}

\section{Existence and uniqueness of the solution}

We shall also need a hypothesis on the initial condition $u_{0}$:

\begin{enumerate}
\item[(\textbf{H.1})] $u_{0}$ is an $\mathcal{G}_{0}$--measurable random
variable s.t. $\sup_{x\in \mathbb{R}^{d}}E|u_{0}(x)|^{p_{0}}<+\infty
$ for some $p_{0}$ large enough in $[2,+\infty [.$
\end{enumerate}

The main result of this section is the following:

\begin{Theo}
\label{teu}Suppose that $b$ and $\sigma $ are Lipschitz functions and assume
that $(\mathbf{H.1})$ and $(\mathbf{H}_{1}^{\alpha })$ are satisfied. Then,
there exists a unique solution $u(t,x)$ to $(\ref{eqfrac})$ such that
\begin{equation}
\sup_{t\in [0,T]}\sup_{x\in \mathbb{R}^{d}}E\left\vert
u(t,x)\right\vert ^{p}<+\infty ,\text{ \ for any }T>0\text{ and
}p\in [2,p_{0}]\text{.}  \label{eq23}
\end{equation}
\end{Theo}

\begin{Rem}
Our results improve those of \cite{DeDo} and \cite{AM} to $d$--dimensional
case and extend the result in \cite{Da} to the fractional setting, both with
$\alpha _{i}$ in $]0,2]\setminus \{1\}$ for $i=1,\ldots ,d$.
\end{Rem}

\noindent \textbf{Proof}. The proof of the existence is done by the Picard
iteration scheme. That is, we define recursively for all $n\geq 1.$
\begin{equation}
\left\{
\begin{array}{l}
u^{0}(t,x)=\displaystyle\int_{\mathbb{R}^{d}}\mathbf{G}_{\alpha ,\delta
}(t,x-y)u_{0}(y)dy, \\
u^{n}(t,x)=u^{0}(t,x)+\displaystyle\int_{0}^{t}\displaystyle\int_{\mathbb{R}%
^{d}}\mathbf{G}_{\alpha ,\delta }(t-s,x-y)\sigma \left( u^{n-1}(s,y)\right)
M(ds,dy) \\
\ \ \ \ \ \ \ \ \ \ \ \ \ \ \ \ +\displaystyle\int_{0}^{t}ds\displaystyle%
\int_{\mathbb{R}^{d}}\mathbf{G}_{\alpha ,\delta }(t-s,x-y)b\left(
u^{n-1}(s,y)\right) dy.%
\end{array}%
\right.  \label{eq24}
\end{equation}%
Fix $p\in [2,p_{0}]$. To prove the theorem we only need the
following facts:

\begin{enumerate}
\item[(i)] For $n\geq 0$, $u^{n}=\{u^{n}(t,x):(t,x)\in [0,T]\times
\mathbb{R}^{d}\}$ is predictable process which satisfies:%
\begin{equation*}
(P_{n})\ \ \sup\limits_{t\in [0,T]}\sup\limits_{x\in \mathbb{R}%
^{d}}E|u^{n}(t,x)|^{2}<+\infty .
\end{equation*}%
This proves that the sequence $u^{n}$ is well defined.

\item[(ii)] For $T>0$ and any $p\in [2,p_{0}]$, $\sup_{n\geq
0}\sup_{t\in [0,T]}\sup_{x\in \mathbb{R}^{d}}E|u^{n}(t,x)|^{p}<+%
\infty $.

\item[(iii)] Set $v_{n}(t)=\sup_{x\in \mathbb{R}^{d}}E\left\vert
u^{n+1}(t,x)-u^{n}(t,x)\right\vert ^{p},\;n\geq 0$. Then
\begin{equation}
v_{n}(t)\leq c\int_{0}^{t}v_{n-1}(s)(J(t-s)+1)ds,  \label{eq25}
\end{equation}%
where
\begin{equation*}
J(t)=\int_{\mathbb{R}^{d}}|\mathcal{F}\mathbf{G}_{\alpha ,\delta }(t,\cdot
)(\xi )|^{2}\mu (d\xi ).
\end{equation*}%
From this we conclude that the sequence $u^{n}(t,x)$ converges in $%
L^{p}(\Omega )$.
\end{enumerate}

\noindent \textit{Statement }(i) is checked by induction on $n$ as a
consequence of the hypotheses $(\mathbf{H.1})$, $(\mathbf{H}_{1}^{\alpha })$
and the properties on the coefficients. Indeed Lemma \ref{holder} applied to
the probability measure $\mathbf{G}_{\alpha ,\delta }(t,x-y)dy$ gives
\begin{eqnarray*}
E|u^{0}(t,x)|^{2} &=&E\left\vert \int_{\mathbb{R}^{d}}\mathbf{G}_{\alpha
,\delta }(t,x-y)u_{0}(y)dy\right\vert ^{2} \\
&\leq &E\int_{\mathbb{R}^{d}}\mathbf{G}_{\alpha ,\delta }(t,x-y)\left\vert
u_{0}(y)\right\vert ^{2}dy \\
&\leq &\sup_{y\in \mathbb{R}^{d}}E\left\vert u_{0}(y)\right\vert ^{2}\leq
\sup_{y\in \mathbb{R}^{d}}E\left\vert u_{0}(y)\right\vert ^{p_{0}}<+\infty .
\end{eqnarray*}%
This proves $(P_{0})$. Now, we assume that the property $(P_{\ell })$ is
true for any integer $\ell \leq n-1$ ($n\geq 2$). We want to check $(P_{n})$%
. Burkholder's inequality, the induction hypothesis the Lipschitz condition
on $\sigma $ and $(\mathbf{H}_{1}^{\alpha })$ imply
\begin{eqnarray}
&&\sup\limits_{t\in [0,T]}\sup\limits_{x\in \mathbb{R}%
^{d}}E\left\vert \int_{0}^{t}\int_{\mathbb{R}^{d}}\mathbf{G}_{\alpha ,\delta
}(t-s,x-y)\sigma \left( u^{n-1}(s,y)\right) M(ds,dy)\right\vert ^{2}  \notag
\\
&\leq &c\sup\limits_{s\in [0,T]}\sup\limits_{y\in \mathbb{R}%
^{d}}E\left( 1+\left\vert u^{n-1}(s,y)\right\vert ^{2}\right)
\int_{0}^{T}ds\int_{\mathbb{R}^{d}}\left\vert \mathcal{F}\mathbf{G}_{\alpha
,\delta }(s,\cdot )(\xi )\right\vert ^{2}\mu (d\xi ).  \notag
\end{eqnarray}%
For the term containing the drift coefficient $b$ we apply Lemma \ref{holder}
with respect to the measure $\mathbf{G}_{\alpha ,\delta }(s,y)dsdy$ to
obtain
\begin{eqnarray*}
&&\sup\limits_{t\in [0,T]}\sup\limits_{x\in \mathbb{R}%
^{d}}E\left\vert \int_{0}^{t}ds\int_{\mathbb{R}^{d}}\mathbf{G}_{\alpha
,\delta }(t-s,x-y)b\left( u^{n-1}(s,y)\right) dy\right\vert ^{2} \\
&\leq &c\sup\limits_{t\in [0,T]}\sup\limits_{x\in \mathbb{R}%
^{d}}E\left( 1+\left\vert u^{n-1}(t,x)\right\vert ^{2}\right) \left(
\int_{0}^{T}ds\int_{\mathbb{R}^{d}}\mathbf{G}_{\alpha ,\delta
}(s,y)dy\right) ^{2}
\end{eqnarray*}%
which is also finite. Hence
\begin{equation*}
\sup\limits_{t\in [0,T]}\sup\limits_{x\in \mathbb{R}%
^{d}}E|u^{n}(t,x)|^{2}\leq c\sup\limits_{t\in
[0,T]}\sup\limits_{x\in \mathbb{R}^{d}}E\left( 1+\left\vert
u^{n-1}(t,x)\right\vert ^{2}\right) ,
\end{equation*}%
we deduce $(P_{n})$.

It is easy to check that $u^{n}$ is predictable.

\noindent \textit{Statement }(ii): Fix $p\in [2,p_{0}]$. We first
prove that for any $n\geq 1$, $t\in [0,T]$
\begin{equation}
\sup_{x\in \mathbb{R}^{d}}E\left\vert u^{n}(t,x)\right\vert ^{p}\leq
c_{1}+c_{2}\int_{0}^{t}\sup_{x\in \mathbb{R}^{d}}E\left\vert
u^{n-1}(s,x)\right\vert ^{p}(J(t-s)+1)ds.  \label{eq32}
\end{equation}%
The arguments are not very far from those used in the proof of (i). Indeed,
we have
\begin{equation}
E\left\vert u^{n}(t,x)\right\vert ^{p}\leq c\left(
C_{0}(t,x)+A_{n}(t,x)+B_{n}(t,x)\right)  \label{eq33}
\end{equation}%
with
\begin{eqnarray*}
C_{0}(t,x) &=&E\left\vert u^{0}(t,x)\right\vert ^{p},\; \\
A_{n}(t,x) &=&E\left\vert \int_{0}^{t}\int_{\mathbb{R}^{d}}\mathbf{G}%
_{\alpha ,\delta }(t-s,x-y)\sigma \left( u^{n-1}(s,y)\right)
M(ds,dy)\right\vert ^{p}, \\
B_{n}(t,x) &=&E\left\vert \int_{0}^{t}ds\int_{\mathbb{R}^{d}}\mathbf{G}%
_{\alpha ,\delta }(t-s,x-y)b\left( u^{n-1}(s,y)\right) dy\right\vert ^{p}.
\end{eqnarray*}%
Jensen's inequality and assumption $(\mathbf{H.1})$ yield
\begin{equation}
\sup\limits_{t\in [0,T]}\sup\limits_{x\in \mathbb{R}%
^{d}}C_{0}(t,x)\leq \sup\limits_{t\in [0,T]}\sup\limits_{x\in \mathbb{%
R}^{d}}\left( E\left\vert u^{0}(t,x)\right\vert ^{p_{0}}\right) ^{\frac{p}{%
p_{0}}}<+\infty .  \label{eq34}
\end{equation}%
Burkholder's inequality and the linear growth condition of $\sigma $ lead to
\begin{equation}
\sup_{x\in \mathbb{R}^{d}}A_{n}(t,x)\leq c\int_{0}^{t}\left( 1+\sup_{y\in
\mathbb{R}^{d}}E\left\vert u^{n-1}(s,y)\right\vert ^{p}\right) J(t-s)ds.
\label{eq35}
\end{equation}%
Moreover, Lemma \ref{holder}, the linear growth condition of $b$ and $(i)$
of Lemma \ref{lgr} imply
\begin{eqnarray}
\sup_{x\in \mathbb{R}^{d}}B_{n}(t,x) &\leq &c\int_{0}^{t}ds\left(
1+\sup_{y\in \mathbb{R}^{d}}E\left\vert u^{n-1}(s,y)\right\vert ^{p}\right)
\int_{\mathbb{R}^{d}}\mathbf{G}_{\alpha ,\delta }(t-s,y)dy  \notag \\
&=&c\int_{0}^{t}\left( 1+\sup_{y\in \mathbb{R}^{d}}E\left\vert
u^{n-1}(s,y)\right\vert ^{p}\right) ds.  \label{eq36}
\end{eqnarray}%
Plugging the estimates (\ref{eq34}) to (\ref{eq36}) into (\ref{eq33}) yields
(\ref{eq32}).

\noindent Finally, the conclusion of part (ii) follows by applying Lemma \ref%
{l2} quoted in the appendix below to the situation: $f_{n}(t)=\sup_{x\in
\mathbb{R}^{d}}E\left\vert u^{n}(t,x)\right\vert ^{p}$, $k_{1}=c_{1}$, $%
k_{2}=0$, $g(s)=c_{2}(J(s)+1)$, with $c_{1},c_{2}$ given in (\ref{eq32}).

\noindent \textit{Statement }(iii): Consider the decomposition
\begin{equation*}
E\left\vert u^{n+1}(t,x)-u^{n}(t,x)\right\vert ^{p}\leq
c(a_{n}(t,x)+b_{n}(t,x)),
\end{equation*}%
with%
\begin{eqnarray*}
a_{n}(t,x) &=&E\left\vert \int_{0}^{t}\int_{\mathbb{R}^{d}}\mathbf{G}%
_{\alpha ,\delta }(t-s,x-y)\left[ \sigma (u^{n}(s,y))-\sigma (u^{n-1}(s,y))%
\right] M(ds,dy)\right\vert ^{p}, \\
b_{n}(t,x) &=&E\left\vert \int_{0}^{t}ds\int_{\mathbb{R}^{d}}\mathbf{G}%
_{\alpha ,\delta }(t-s,x-y)\left[ b(u^{n}(s,y))-b(u^{n-1}(s,y))\right]
\right\vert ^{p}.
\end{eqnarray*}%
Again Burkholder's inequality, Lemma \ref{holder} and the Lipschitz property
of $b$ and $\sigma $ imply that
\begin{equation*}
E\left\vert u^{n+1}(t,x)-u^{n}(t,x)\right\vert ^{p}\leq c\int_{0}^{t}\left(
\sup_{y\in \mathbb{R}^{d}}E\left\vert u^{n}(s,y)-u^{n-1}(s,y)\right\vert
^{p}\right) \left( J(t-s)+1\right) ds.
\end{equation*}%
This yields (\ref{eq25}). We finish the proof by applying Lemma \ref{l2} in
the particular case: $f_{n}(t)=v_{n}(t)$, $k_{1}=k_{2}=0$, $g(s)=c(J(s)+1)$,
with $c$ given in (\ref{eq25}). Notice that the results proved in part (ii)
show that $v:=\sup_{0\leq s\leq T}f_{0}(s)$ is finite. Hence the series $%
\sum_{n\geq 0}(v_{n}(t))^{1/p}$ converges for any $p\in [2,p_{0}]$.
We then conclude that $\{u_{n}(t,x):(t,x)\in [0,T]\times \mathbb{R}%
^{d}\}$ converges uniformly in $L^{p}(\Omega )$ to a limit $%
u=\{u(t,x):(t,x)\in [0,T]\times \mathbb{R}^{d}\}$. It is not
difficult to check that $u$ satisfies conditions of the Definition
\ref{d1} and therefore the existence is completely proved.\hfill
$\square $

Let us now give a short proof for the uniqueness.

Let $u_{1}$ and $u_{2}$ be two solutions to (\ref{eqfrac})%
\begin{equation*}
{F(t,x):=E}\left\vert {u_{1}(t,x)-u_{2}(t,x)}\right\vert {^{2}}
\end{equation*}%
and ${H(t)=\sup_{x\in \mathbb{R}^{d}}F(t,x)}$. Then the isometry property
for stochastic integrals, Lemma \ref{holder} and the Lipschitz condition on $%
b$ and $\sigma $ yield
\begin{equation*}
H(t)\leq c\int_{0}^{t}H(s)(J(t-s)+1)ds.
\end{equation*}%
By iterating this formula and using Lemma \ref{l2} we deduce that ${H\equiv 0%
}$, hence, $u_{1}(t,x)=u_{2}(t,x)$ $t,x$ a.e.\hfill $\square $

\section{Path regularity of the solution}

In this section we analyze the path H\"{o}lder regularity of $u(t,x)$. The
next Theorem extends and improves similar results known for stochastic heat
equation, corresponding to the case $\alpha _{i}=2$, $\delta _{i}=0$, $%
i=1,\ldots ,d$. (see for instance \cite{SS}). Contrary to the factorization
method usually used in high dimension (see e.g. \cite{Sa} and \cite{DaZa})
we use a direct method to prove our regularity results; in which the Fourier
transform and the representation of the Green kernel play a crucial role.

Let $\gamma =(\gamma _{1},\gamma _{2})$ such that $\gamma _{1},\gamma _{2}>0$
and let $K$ be a compact set of $\mathbb{R}^{d}$. We denote by $\mathcal{C}%
^{\gamma }\left( [0,T]\times K;\mathbb{R}\right) $ the set of $\gamma $--H%
\"{o}lder continuous functions equipped with the norm defined by:
\begin{equation*}
\left\Vert f\right\Vert _{\gamma ,T,K}=\sup_{(t,x)\in [0,T]\times
K}\left\vert f(t,x)\right\vert +\sup_{s\neq t\in [0,T]}\sup_{x\neq
y\in K}\dfrac{\left\vert f(t,x)-f(s,y)\right\vert }{\left\vert
t-s\right\vert ^{\gamma _{1}}+\left\Vert x-y\right\Vert ^{\gamma
_{2}}}.
\end{equation*}

\begin{enumerate}
\item[(\textbf{H.2})] There exists $\rho \in ]0,1[$ s.t. for all $z\in K$
(compact subset of $\mathbb{R}^{d}$)
\begin{equation*}
\sup_{y\in \mathbb{R}^{d}}E|u_{0}(y+z)-u_{0}(y)|^{p_{0}}\leq c|z|^{\rho
p_{0}}
\end{equation*}%
for some $p_{0}$ large enough in $[2,\tfrac{\alpha _{0}}{\rho }[$, where $%
\alpha _{0}=\min_{1\leq i\leq d}\{\alpha _{i}\}$.
\end{enumerate}

The main result of this section is

\begin{Theo}
\label{thr}Suppose that $b$ and $\sigma $ are Lipschitz functions. Assume
that $(\mathbf{H.1})$, $(\mathbf{H.2})$ and $(\mathbf{H}_{\eta }^{\alpha })$
hold for $\eta \in ]0,1[$ and let $u$ be a solution to equation $(\ref%
{eqfrac})$. Then $u$ belongs a.s. to the H\"{o}lder space $\mathcal{C}%
^{\gamma }\left( [0,T]\times K;\mathbb{R}\right) $ for $0<\gamma _{1}<\min
\{\sum_{i=1}^{d}(\rho /\alpha _{i}),(1-\eta )/2\}$, $0<\gamma _{2}<\min
\{\rho ,\alpha _{0}(1-\eta )/2,1/2\}$ and $K$ compact subset of $\mathbb{R}%
^{d}$. Moreover $E\left\Vert u\right\Vert _{\gamma ,T,K}^{p}<+\infty
$ for any $p\in [2,p_{0}]$.
\end{Theo}

In what follows, we establish some technical and useful results that will be
needed in the proof of the regularity of the solution.

Let $u_{0}$ be a given random fields satisfying $\mathbf{(H.1)}$. Set
\begin{equation*}
u^{0}(t,x)={\int_{\mathbb{R}^{d}}}\mathbf{G}_{\alpha ,\delta
}(t,x-y)u_{0}(y)dy.
\end{equation*}

\begin{Pro}
\label{pinitial}Suppose that $u_{0}$ satisfies $\mathbf{(H.2})$ then

\begin{description}
\item[$(i)$] For each $x\in \mathbb{R}^{d}$ a.s. the mapping $t\longmapsto
u^{0}(t,x)$ is $\gamma _{1}$--H\"{o}lder continuous for $0<\gamma _{1}<\min
\{\sum_{i=1}^{d}(\rho /\alpha _{i}),1\}.$

\item[$(ii)$] For each $t\in [0,T]$ a.s. the mapping $x\longmapsto
u^{0}(t,x)$ is $\gamma _{2}$--H\"{o}lder continuous for $0<\gamma _{2}<\rho $%
.
\end{description}
\end{Pro}

\noindent \textbf{Proof. }\newline
\noindent \textit{Proof of (i)}\textbf{\ }Using the semigroup property of
the Green kernel $\mathbf{G}_{\alpha ,\delta }$,
\begin{equation*}
\mathbf{G}_{\alpha ,\delta }(t+h,x-y)=\int_{\mathbb{R}^{d}}\mathbf{G}%
_{\alpha ,\delta }(t,x-y-z)\mathbf{G}_{\alpha ,\delta }(h,z)dz,
\end{equation*}%
we have
\begin{eqnarray*}
u^{0}(t+h,x)-u^{0}(t,x)={\int_{\mathbb{R}^{d}}}\left[ \mathbf{G}_{\alpha
,\delta }(t+h,x-y)-\mathbf{G}_{\alpha ,\delta }(t,x-y)\right] u_{0}(y)dy &&
\\
={\int_{\mathbb{R}^{d}}}\mathbf{G}_{\alpha ,\delta }(h,z)\left( {\int_{%
\mathbb{R}^{d}}}\mathbf{G}_{\alpha ,\delta }(t,x-y)\left[ u_{0}(y-z)-u_{0}(y)%
\right] dy\right) dz. &&
\end{eqnarray*}%
By concavity of the of the mapping $x\longmapsto x^{p/p_{0}}$ (since
$p\leq p_{0})$, Lemma \ref{holder}, the assumption $\mathbf{(H.2)}$
and Fubini's theorem we obtain for $p\in [2,p_{0}]$
\begin{eqnarray*}
E\left\vert u^{0}(t+h,x)-u^{0}(t,x)\right\vert ^{p} &\leq &\left( {\int_{%
\mathbb{R}^{d}}}\sup_{y\in \mathbb{R}^{d}}E|u_{0}(y-z)-u_{0}(y)|^{p_{0}}%
\mathbf{G}_{\alpha ,\delta }(h,z)dz\right) ^{\frac{p}{p_{0}}} \\
&\leq &c\left( {\int_{\mathbb{R}^{d}}}\mathbf{G}_{\alpha ,\delta
}(h,z)|z|^{\rho p_{0}}dz\right) ^{\frac{p}{p_{0}}}.
\end{eqnarray*}%
Now, due to the scaling property of $\mathbf{G}_{\alpha ,\delta }$
\begin{equation*}
\int_{\mathbb{R}^{d}}\mathbf{G}_{\alpha ,\delta }(h,z)|z|^{\rho
p_{0}}dz=\int_{\mathbb{R}^{d}}h^{-\sum_{i=1}^{d}\frac{1}{\alpha _{i}}}%
\mathbf{G}_{\alpha ,\delta }(1,h^{-\sum_{i=1}^{d}\frac{1}{\alpha _{i}}%
}z)|z|^{\rho p_{0}}dz
\end{equation*}%
\begin{equation*}
=h^{\sum_{i=1}^{d}\frac{\rho p_{0}}{\alpha _{i}}}\int_{\mathbb{R}^{d}}%
\mathbf{G}_{\alpha ,\delta }(1,y)|y|^{\rho p_{0}}dy.
\end{equation*}%
Using property $(iv)$ of Lemma \ref{lgr} we obtain that
\begin{eqnarray*}
\int_{\mathbb{R}^{d}}\mathbf{G}_{\alpha ,\delta }(1,y)|y|^{\rho p_{0}}dy
&\leq &c\sum_{i=1}^{d}\int_{\mathbb{R}}\mathbf{G}_{\alpha ,\delta
}(1,y_{i})|y_{i}|^{\rho p_{0}}dy \\
&\leq &c\sum_{i=1}^{d}\int_{\mathbb{R}}\frac{|y|^{\rho p_{0}}}{%
1+|y|^{1+\alpha _{i}}}dy \\
&\leq &c\int_{0}^{+\infty }\frac{x^{\rho p_{0}}}{1+x^{1+\alpha _{0}}}dy.
\end{eqnarray*}%
The last integral is convergent for $\rho p_{0}<\alpha _{0}$. Therefore we
have proved the assertion (i).

\noindent \textit{Proof of (ii) }A change of variable gives immediately
\begin{eqnarray*}
u^{0}(t,x+z)-u^{0}(t,x) &=&{\int_{\mathbb{R}^{d}}}\left[ \mathbf{G}_{\alpha
,\delta }(t,x+z-y)-\mathbf{G}_{\alpha ,\delta }(t,x-y)\right] u_{0}(y)dy \\
&=&{\int_{\mathbb{R}^{d}}}\mathbf{G}_{\alpha ,\delta }(t,x-y)\left[
u_{0}(y+z)-u_{0}(y)\right] dy.
\end{eqnarray*}%
Applying again concavity of $x\longmapsto x^{p/p_{0}}$, Lemma
\ref{holder} for the probability measure $\mathbf{G}_{\alpha ,\delta
}(t,x-y)dy$, the assumption $\mathbf{(H.2)}$ and Fubini's theorem,
we obtain for all $p\in [2,p_{0}]$
\begin{eqnarray*}
&&E\left\vert u^{0}(t,x+z)-u^{0}(t,x)\right\vert ^{p}\leq \\
&\leq &c\left( {\int_{\mathbb{R}^{d}}}\sup_{y\in \mathbb{R}%
^{d}}E(|u_{0}(y+z)-u_{0}(y)|^{p_{0}})\mathbf{G}_{\alpha ,\delta
}(t,x-y)dy\right) ^{\frac{p}{p_{0}}} \\
&\leq &c\left( {\int_{\mathbb{R}^{d}}}\mathbf{G}_{\alpha ,\delta
}(t,x-y)|z|^{\rho p_{0}}dy\right) ^{\frac{p}{p_{0}}}=c|z|^{\rho p}.
\end{eqnarray*}%
We complete the proof of\ the Lemma by Kolmogorov's criterion.\hfill $%
\square $

The next proposition studies the H\"{o}lder regularity of stochastic
integrals with respect to the martingale measure $M$. For a given
predictable random field $V$ we set
\begin{equation*}
U(t,x)=\int_{0}^{t}{\int_{\mathbb{R}^{d}}}\mathbf{G}_{\alpha ,\delta
}(t-s,x-y)V(s,y)M(ds,dy).
\end{equation*}

\begin{Pro}
\label{psi}Assume that $\sup_{0\leq t\leq T}\sup_{x\in \mathbb{R}%
^{d}}E(|V(t,x)|^{p})$ is finite for some $p$ large enough. Then under $%
\mathbf{(H}_{\eta }^{\alpha })$ we have

\begin{description}
\item[$(i)$] For each $x\in \mathbb{R}^{d}$ a.s. the mapping $t\longmapsto
U(t,x)$ is $\gamma _{1}$--H\"{o}lder continuous for $0<\gamma _{1}<(1-\eta
)/2.$

\item[$(ii)$] For each $t\in [0,T]$ a.s. the mapping $x\longmapsto
U(t,x)$ is $\gamma _{2}$--H\"{o}lder continuous for $0<\gamma _{2}<\min
\{\alpha _{0}(1-\eta )/2,1/2\}$.
\end{description}
\end{Pro}

\noindent \textbf{Proof. }\newline
\noindent \textit{Proof of (i). }We have%
\begin{eqnarray*}
U(t+h,x)-U(t,x)=\int_{0}^{t+h}{\int_{\mathbb{R}^{d}}}\mathbf{G}_{\alpha
,\delta }(t+h-s,x-y)V(s,y)M(ds,dy) && \\
-\int_{0}^{t}{\int_{\mathbb{R}^{d}}}\mathbf{G}_{\alpha ,\delta
}(t-s,x-y)V(s,y)M(ds,dy) && \\
=\int_{0}^{t}{\int_{\mathbb{R}^{d}}}\left[ \mathbf{G}_{\alpha ,\delta
}(t+h-s,x-y)-\mathbf{G}_{\alpha ,\delta }(t-s,x-y)\right] V(s,y)M(ds,dy) &&
\\
+\int_{t}^{t+h}{\int_{\mathbb{R}^{d}}}\mathbf{G}_{\alpha ,\delta
}(t+h-s,x-y)V(s,y)M(ds,dy). &&
\end{eqnarray*}%
For $p\geq 2$, the Burkholder inequality yields
\begin{eqnarray*}
&&E\left\vert U(t+h,x)-U(t,x)\right\vert ^{p} \\
&\leq &c\sup_{(s,y)\in [0,T]\times \mathbb{R}^{d}}E\left\vert
V(s,y)\right\vert ^{p}\left[ \left( I_{1}(h,x)\right)
^{\frac{p}{2}}+\left( I_{2}(h,x)\right) ^{\frac{p}{2}}\right] ,
\end{eqnarray*}%
where
\begin{equation*}
I_{1}(h,x)=\int_{0}^{T}ds\int_{\mathbb{R}^{d}}\left\vert \mathcal{F}\left[
\mathbf{G}_{\alpha ,\delta }(h+s,x-\cdot )-\mathbf{G}_{\alpha ,\delta
}(s,x-\cdot )\right] (\xi )\right\vert ^{2}\mu (d\xi )
\end{equation*}%
and
\begin{equation*}
I_{2}(h,x)=\int_{0}^{h}ds\int_{\mathbb{R}^{d}}\left\vert \mathcal{F}\mathbf{G%
}_{\alpha ,\delta }(s,x-\cdot )(\xi )\right\vert ^{2}\mu (d\xi ).
\end{equation*}%
\textit{Estimation of}\textbf{\ }$I_{1}(h,x)$. Set
\begin{equation}
\psi _{\alpha ,\xi }(t):=\exp \{-t\sum_{i=1}^{d}|\xi _{i}|^{\alpha _{i}}\exp
(-\imath \delta _{i}\frac{\pi }{2}\mathrm{sgn}(\xi _{i}))\}.  \label{eqpsi}
\end{equation}%
By the definition of the Green function $\mathbf{G}_{\alpha ,\delta }$ we
can write
\begin{eqnarray*}
I_{1}(h,x) &=&\int_{0}^{T}dr\int_{\mathbb{R}^{d}}\left\vert \psi _{\alpha
,\xi }(r)\psi _{\alpha ,\xi }(h)-\psi _{\alpha ,\xi }(r)\right\vert ^{2}\mu
(d\xi ) \\
&=&\int_{0}^{T}dr\int_{\mathbb{R}^{d}}\left\vert \psi _{\alpha ,\xi
}(r)\right\vert ^{2}\left\vert \psi _{\alpha ,\xi }(h)-1\right\vert ^{2}\mu
(d\xi ).
\end{eqnarray*}%
From the proof of Lemma \ref{lemmu} we know that $|\psi _{\alpha ,\xi
}(r)|^{2}\leq \exp (-2r\kappa S_{\alpha }(\xi ))$. The mean value theorem to
the function $\psi _{\alpha ,\xi }(\cdot )$ between $0$ and $h$ implies that
\begin{equation*}
\left\vert \psi _{\alpha ,\xi }(h)-1\right\vert \leq hS_{\alpha }(\xi ).
\end{equation*}%
\text{Moreover, for any }$\beta \in ]0,1[$, we have $|\psi _{\alpha ,\xi
}(h)-1|\leq 2^{1-\beta }|\psi _{\alpha ,\xi }(h)-1|^{\beta }.$\newline
Therefore
\begin{equation*}
I_{1}(h,x)\leq c_{\beta }h^{2\beta }\int_{0}^{T}dr\int_{\mathbb{R}%
^{d}}(S_{\alpha }(\xi ))^{2\beta }\exp \left( -2\kappa rS_{\alpha }(\xi
)\right) \mu (d\xi ).
\end{equation*}%
Hence, under $\mathbf{(H}_{\eta }^{\alpha })$ Lemma \ref{lap1} implies that
the right hand side of the above inequality is finite for any $\beta $ in $%
]0,(1-\eta )/2[$.

\noindent \textit{Estimation of}\textbf{\ }$I_{2}(h,x)$. Fubini's theorem
and a change of variable lead to
\begin{eqnarray*}
I_{2}(h,x) &=&\int_{\mathbb{R}^{d}}\left[ \int_{0}^{h}ds\left\vert \mathcal{F%
}\mathbf{G}_{\alpha ,\delta }(s,x-\cdot )(\xi )\right\vert ^{2}\right] \mu
(d\xi \\
&=&\int_{\mathbb{R}^{d}}\left[ \int_{0}^{h}\left\vert \psi _{\alpha ,\xi
}(s)\right\vert ^{2}ds\right] \mu (d\xi ) \\
&\leq &\int_{\mathbb{R}^{d}}\left[ \int_{0}^{h}\exp \left( -2\kappa
rS_{\alpha }(\xi )\right) dr\right] \mu (d\xi ).
\end{eqnarray*}%
By Lemma \ref{lap2} the last term is bounded by $c(h+h^{\beta })$ for any $%
\beta \in ]0,1-\eta [$.

\noindent \textit{Proof of (ii)}. Let $x\in \mathbb{R}^{d}$ and $z$ belongs
to a compact subset $K$ of $\mathbb{R}^{d}$
\begin{eqnarray*}
U(t,x+z)-U(t,x) &=&\int_{0}^{t}{\int_{\mathbb{R}^{d}}}\left[ \mathbf{G}%
_{\alpha ,\delta }(t-s,x+z-y)-\mathbf{G}_{\alpha ,\delta }(t-s,x-y)\right] \\
&&\ \ \ \ \ \ \ \ \ \ \ \times V(s,y)M(ds,dy).
\end{eqnarray*}%
For $p\geq 2$, Burkholder's inequality yields
\begin{equation*}
E\left\vert U(t,x+z)-U(t,x)\right\vert ^{p}\leq c\sup_{(s,y)\in
[0,T]\times \mathbb{R}^{d}}E\left\vert V(s,y)\right\vert ^{p}\left(
J(x,z)\right) ^{\frac{p}{2}},
\end{equation*}%
where
\begin{equation*}
J(x,z)=\int_{0}^{T}ds\int_{\mathbb{R}^{d}}\left\vert \mathcal{F}\left[
\mathbf{G}_{\alpha ,\delta }(s,x+z-\cdot )-\mathbf{G}_{\alpha ,\delta
}(s,x-\cdot )\right] (\xi )\right\vert ^{2}\mu (d\xi )
\end{equation*}%
The property $\mathcal{F}\left( f\right) (\xi +a)=\mathcal{F}\left(
e^{-2i\pi \left\langle a,\cdot \right\rangle }f\right) (\xi )$ together with
Lemma \ref{lap3} imply
\begin{eqnarray*}
J(x,z) &=&\int_{0}^{T}ds\int_{\mathbb{R}^{d}}\left\vert \psi _{\alpha ,\xi
}(s)\left[ \exp \left( 2\imath \pi \left\langle \xi ,z\right\rangle \right)
-1\right] \right\vert ^{2}\mu (d\xi ) \\
&\leq &4\int_{\mathbb{R}^{d}}\int_{0}^{T}dr\exp \left( -2\kappa rS_{\alpha
}(\xi )\right) \left[ 1-\cos \left( 2\pi \left\langle \xi ,z\right\rangle
\right) \right] \mu (d\xi ) \\
&\leq &c(\left\Vert z\right\Vert _{\infty }+\left\Vert z\right\Vert _{\infty
}^{2\beta }).
\end{eqnarray*}%
for all $0<\beta <\min \{\alpha _{0}(1-\eta )/2,1/2\}$%
.\hfill $\square $

\bigskip

\noindent \textbf{Proof of Theorem \ref{thr}.} Let $u$ be a solution
to equation (\ref{eqfrac}). For any $x\in \mathbb{R}^{d}$, $z\in K$
and $t\in [0,T]$, and for $p\in [2,p_{0}]$,
\newpage
\begin{eqnarray*}
&&E|u(t+h,x+z)-u(t,x)|^{p}\leq cE|u^{0}(t+h,x+z)-u^{0}(t,x)|^{p} \\
&&+cE\left( \left\vert \int_{0}^{t+h}\int_{\mathbb{R}^{d}}\mathbf{G}_{\alpha
,\delta }(t+h-s,x+z-y)\sigma \left( u(s,y)\right) M(ds,dy)\right. \right. \\
&&\ \ \ \ \ \ \ \ \ \ \ \ \ \ \left. \left. -\int_{0}^{t}\int_{\mathbb{R}%
^{d}}\mathbf{G}_{\alpha ,\delta }(t-s,x-y)\sigma \left( u(s,y)\right)
M(ds,dy)\right\vert \right) ^{p} \\
&&+cE\left( \left\vert \int_{0}^{t+h}\int_{\mathbb{R}^{d}}\mathbf{G}_{\alpha
,\delta }(t+h-s,x+z-y)b\left( u(s,y)\right) dsdy\right. \right. \\
&&\ \ \ \ \ \ \ \ \ \ \ \ \ \ \left. \left. -\int_{0}^{t}\int_{\mathbb{R}%
^{d}}\mathbf{G}_{\alpha ,\delta }(t-s,x-y)b\left( u(s,y)\right)
dsdy\right\vert \right) ^{p} \\
&=&A_{1}+A_{2}+A_{3}.
\end{eqnarray*}%
The terms $A_{1}$, $A_{2}$ are estimated by Proposition \ref{pinitial} and %
\ref{psi}. More precisely we get
\begin{equation*}
A_{1}+A_{2}\leq c\left( h^{\sum_{i=1}^{d}\frac{\rho p}{\alpha _{i}}%
}+h^{\beta _{1}p}+\left\Vert z\right\Vert _{\infty }^{\rho p}+\left\Vert
z\right\Vert _{\infty }^{\beta _{2}p}\right)
\end{equation*}%
for $\beta _{1}<\frac{1-\eta }{2}$ and $\beta _{2}<\min (\rho ,\frac{\alpha
_{0}(1-\eta )}{2},\frac{1}{2})$.\newline
Let us now give the estimation of $A_{3}$. After a change of variable $%
A_{3}=|B|^{p}$, where
\begin{eqnarray*}
B &=&\int_{0}^{h}\int_{\mathbb{R}^{d}}\mathbf{G}_{\alpha ,\delta
}(h+t-s,x+z-y)b\left( u(s,y)\right) dsdy \\
&&+\int_{0}^{t}\int_{\mathbb{R}^{d}}\mathbf{G}_{\alpha ,\delta }(t-s,x-y)%
\left[ b\left( u(s+h,y+z)\right) -b\left( u(s,y)\right) \right] dsdy.
\end{eqnarray*}%
By Lemma \ref{holder} to the measure $\mathbf{G}_{\alpha ,\delta }(s,y)dsdy$%
, Lipschitz and linear growth condition of $b$, we deduce
\begin{eqnarray*}
A_{3} &\leq &\int_{0}^{h}\int_{\mathbb{R}^{d}}\mathbf{G}_{\alpha ,\delta
}(t+h-s,x+z-y)E\left\vert b\left( u(s,y)\right) \right\vert ^{p}dsdy \\
&&+\int_{0}^{t}\int_{\mathbb{R}^{d}}\mathbf{G}_{\alpha ,\delta
}(t-s,x+z-y)E\left\vert b\left( u(s+h,y+z)\right) -b\left( u(s,y)\right)
\right\vert ^{p}dsdy \\
&\leq &c_{1}\cdot h+c_{2}\int_{0}^{t}\sup_{y\in \mathbb{R}^{d}}E\left\vert
u(s+h,y+z)-u(s,y)\right\vert ^{p}ds,
\end{eqnarray*}%
where we have used assertion (i) of Lemma \ref{lgr}.

Put $\varphi (s,h,z)=\sup_{y\in \mathbb{R}^{d}}E\left\vert
u(s+h,y+z)-u(s,y)\right\vert ^{p}$. Hence,
\begin{equation*}
\varphi (t,h,z)\leq c_{3}(h+h^{\sum_{i=1}^{d}\frac{\rho p}{\alpha _{i}}%
}+h^{\beta _{1}p}+\left\Vert z\right\Vert _{\infty }^{\rho p}+\left\Vert
z\right\Vert _{\infty }^{\beta _{2}p})+c_{2}\int_{0}^{t}\varphi (s,h,z)ds.
\end{equation*}%
Therefore by Gronwall's lemma
\begin{eqnarray*}
\sup_{0\leq t\leq T}\sup_{x\in \mathbb{R}^{d}}E\left\vert
u(t+h,x+z)-u(t,x)\right\vert ^{p}\leq && \\
\leq c(h+h^{\sum_{i=1}^{d}\frac{\rho p}{\alpha _{i}}}+h^{\beta
_{1}p}+\left\Vert z\right\Vert _{\infty }^{\rho p}+\left\Vert z\right\Vert
_{\infty }^{\beta _{2}p}). &&
\end{eqnarray*}%
The conclusion of Theorem \ref{thr} is a consequence of the Kolmogorov
continuity criterium.\hfill $\square $

\begin{Rem}
Note that our results of regularity of the solution generalizes those
obtained in \cite{SS}.
\end{Rem}

\section{Smoothness of the law}

We prove that, under non--degeneracy condition on the diffusion coefficient $%
\sigma $, the law of $u(t,x)$ (solution of $Eq_{\delta }^{\alpha
}(d,b,\sigma )$), has a smooth density for fixed $t>0$, $x\in \mathbb{R}^{d}$%
.

\noindent Since the noise $F$ has a space correlation, the setting of the
corresponding stochastic calculus of variations is that used in \cite{MS}
(see also \cite{MMS}).

\noindent Let $\mathcal{E}$ be the inner--product space consisting of
functions $\varphi :\mathbb{R}^{d}\rightarrow \mathbb{R}$ such that
\begin{equation*}
\int_{\mathbb{R}^{d}}\Gamma (dx)(\varphi \ast \tilde{\varphi})(x)<+\infty ,
\end{equation*}%
where $\tilde{\varphi}(x)=\varphi (-x).$ This space is endowed with the
inner product
\begin{equation*}
\langle \varphi ,\psi \rangle _{\mathcal{E}}=\int_{\mathbb{R}^{d}}\Gamma
(dx)(\varphi \ast \tilde{\psi})(x)=\int_{\mathbb{R}^{d}}\Gamma (dx)\int_{%
\mathbb{R}^{d}}dy\varphi (x-y)\psi (-y).
\end{equation*}%
Let $\mathcal{H}$ be the completion of $\mathcal{E}$ with respect the norm
derived from $\langle \cdot ,\cdot \rangle _{\mathcal{E}}$. Set $\mathcal{H}%
_{T}=L^{2}([0,T];\mathcal{H})$, this space is a real separable Hilbert space
isomorphic to the reproducing kernel Hilbert space of the centred Gaussian
noise $F$ which can be identified with a Gaussian process $\{W(h):h\in
\mathcal{H}_{T}\}$ as follows. For any complete orthonormal system $%
\{e_{k}:k\in \mathbb{N}\}\subset \mathcal{E}$ of the Hilbert space $\mathcal{%
H}$, define
\begin{equation*}
W_{k}(t)=\int_{0}^{t}\int_{\mathbb{R}^{d}}e_{k}(x)F(ds,dx),\ k\in \mathbb{N}%
,\ \ t\in [0,T],
\end{equation*}%
where the integral must be understood in Dalang's sense. The process $%
\{W_{k}(t):t\in [0,T],k\in \mathbb{N}\}$ is a sequence of
independent standard Brownian motions, such that for any predictable
process $X$
\begin{equation*}
\int_{0}^{T}\int_{\mathbb{R}^{d}}X(s,x)F(ds,dx)=\sum_{k=0}^{\infty
}\int_{0}^{T}\left\langle X(s,\cdot ),e_{k}(\cdot )\right\rangle _{\mathcal{H%
}}dW_{k}(s).
\end{equation*}%
In particular, for any $\varphi \in \mathcal{S}(\mathbb{R}^{d+1})$
\begin{equation*}
F(\varphi )=\sum_{k=0}^{\infty }\int_{0}^{T}\left\langle \varphi (s,\cdot
),e_{k}(\cdot )\right\rangle _{\mathcal{H}}dW_{k}(s).
\end{equation*}%
For $h\in \mathcal{H}_{T}$, set
\begin{equation*}
W(h)=\sum_{k=1}^{\infty }\int_{0}^{T}\left\langle h(s),e_{k}\right\rangle _{%
\mathcal{H}}dW_{k}(s).
\end{equation*}%
Thus, we can use the differential Malliavin calculus based on $\{W(h):h\in
\mathcal{H}_{T}\}$ (for more details see \cite{N}). Recall that the Sobolev
spaces $\mathbb{D}^{k,p}$ are defined by means of iterations of the
derivative operator $D$. For a random variable $X$, $D^{k}X$ defines a $%
\mathcal{H}_{T}^{\otimes k}$--valued random variable if it exists. For $h\in
\mathcal{H}_{T}$, set $D_{h}X=\left\langle DX,h\right\rangle _{\mathcal{H}%
_{T}}$ and for $r\in [0,T]$, $D_{r}X$ defines an element of $\mathcal{%
H}$, which is denoted by $D_{r,\cdot }X.$ Then for any $h\in \mathcal{H}_{T}$%
\begin{equation*}
D_{h}X=\int_{0}^{T}\left\langle D_{r,\cdot }X,h(r)\right\rangle _{\mathcal{H}%
}dr.
\end{equation*}%
We write $D_{r,\varphi }X=\left\langle D_{r,\cdot }X,\varphi \right\rangle _{%
\mathcal{H}}$ for $\varphi \in \mathcal{H}$.

The main result in this section is

\begin{Theo}
\label{theo2}Fix $t>0$ and $x\in \mathbb{R}^{d}$. Assume that the
coefficients $\sigma $ and $b$ are $\mathcal{C}^{\infty }$ with bounded
derivatives of any order and in addition there exists $a>0$ such that $%
\sigma (z)\geq a$ for any $z\in \mathbb{R}$. Then under $(\mathbf{H}_{\eta
}^{\alpha })$ with $\eta \in (0,\frac{1}{2})$ the law of $u(t,x)$, solution
of $(\ref{eqfrac})$, has a density which is infinitely differentiable.
\end{Theo}

\begin{Rem}
If $d=1$ and $F$ is a space--time white noise, we obtain existence and
smoothness of the density of the law for solution to equations studied in
\cite{DeDo} and \cite{AM}. In the case where $\alpha _{i}=2$, $\delta _{i}=0$%
, $i=1,\ldots ,d$, we obtain the result obtained in \cite{MMS}.
\end{Rem}

Classical results on the existence and smoothness of the density, using the
approach of Malliavin calculus, are based on the next proposition.

\begin{Pro}
\label{outil}\cite{RoSa} Let $F=(F^{1},\ldots ,F^{m})$, $m\geq 1$, be a
random vector satisfying the following conditions:\newline
$(i)$ $F^{i}$ belongs to $\mathbb{D}^{\infty }=\cap _{p\geq 1}\cap _{k\geq 1}%
\mathbb{D}^{k,p}$ for all $i=1,\ldots ,m.$\newline
$(ii)$ The Malliavin matrix $\gamma _{F}=(\langle F^{i},F^{j}\rangle _{%
\mathcal{H}})_{1\leq i,j\leq m}$ satisfies
\begin{equation*}
(\det \gamma _{F})^{-1}\in \bigcap_{p\geq 1}L^{p}(\Omega ).
\end{equation*}%
Then, $F$ has an infinitely differentiable density with respect to the
Lebesgue measure on $\mathbb{R}^{m}$.
\end{Pro}

According to Theorem 2.1 of \cite{MMS} the proof of Theorem \ref{theo2}, is
a consequence of the following Lemma.

\begin{Lem}
Let $t\in ]0,T]$.

\begin{description}
\item[$(i)$] There exists a positive constant $c_{1}$ such that for any $%
\rho \in [0,t\wedge 1]$,
\begin{equation*}
\int_{0}^{\rho }ds\int_{\mathbb{R}^{d}}\mu (d\xi )|\mathcal{F}\mathbf{G}%
_{\alpha ,\delta }(s,\cdot )(\xi )|^{2}\geq c_{1}\rho ^{\theta _{1}}
\end{equation*}%
for any $\theta _{1}\geq 1$.

\item[$(ii)$] Assume that $(\mathbf{H}_{\eta }^{\alpha })$ holds for $\eta
\in ]0,1[$. Then there exists a positive constant $c_{2}$ such that for any $%
\rho \in [0,t\wedge 1]$,
\begin{equation*}
\int_{0}^{\rho }ds\int_{\mathbb{R}^{d}}\mu (d\xi )|\mathcal{F}\mathbf{G}%
_{\alpha ,\delta }(s,\cdot )(\xi )|^{2}\leq c_{2}\rho ^{\theta _{2}},
\end{equation*}%
for any $\theta _{2}\in ]0,1-\eta ]$.
\end{description}
\end{Lem}

\noindent \textbf{Proof.}\newline
(i) Using\ the lower bound of (\ref{eqfourier}), we have,
\begin{eqnarray*}
\int_{0}^{\rho }ds\int_{\mathbb{R}^{d}}\mu (d\xi )|\mathcal{F}\mathbf{G}%
_{\alpha ,\delta }(s,\cdot )(\xi )|^{2} &\geq &\int_{\mathbb{R}^{d}}\frac{%
\rho \mu (d\xi )}{1+2\rho S_{\alpha }(\xi )} \\
&\geq &\rho \int_{\mathbb{R}^{d}}\frac{\mu (d\xi )}{1+S_{\alpha }(\xi )}\geq
c_{1}\rho ^{\theta _{1}},
\end{eqnarray*}%
for any $\theta _{1}\geq 1$.\newline (ii) Using the upper bound of
(\ref{eqfourier}), we obtain for any $\rho \in [0,t\wedge 1]$,
\begin{eqnarray}
\int_{0}^{\rho }ds\int_{\mathbb{R}^{d}}\mu (d\xi )|\mathcal{F}\mathbf{G}%
_{\alpha ,\delta }(s,\cdot )(\xi )|^{2} &\leq &2\rho \int_{\mathbb{R}^{d}}%
\frac{\mu (d\xi )}{1+2\rho \kappa \sum_{i=1}^{d}|\xi |^{\alpha _{i}}}  \notag
\\
&\leq &c\rho ^{1-\eta }\int_{\mathbb{R}^{d}}\frac{\mu (d\xi )}{%
(1+\sum_{i=1}^{d}|\xi |^{\alpha _{i}})^{\eta }}  \notag
\end{eqnarray}%
and we complete the proof of (ii) by the assumption $(\mathbf{H}_{\eta
}^{\alpha })$.\hfill $\square $

\appendix

\section{Appendix}

In this last section, we collect some of the technical results which have
been used the in previous sections.

\begin{Lem}
\label{holder}Let $f$, $h$ be two functions defined on $\mathbb{R}^{d}$ and $%
\mu $ a positive measure such that $f\cdot h\in L^{1}(\mu )$. Then, for all $%
q>1$, we have:
\begin{equation*}
\left\vert \int f\cdot |h|d\mu \right\vert ^{q}\leq \left( \int |f|^{q}\cdot
|h|d\mu \right) \left( \int |h|d\mu \right) ^{q-1}.
\end{equation*}
\end{Lem}

\noindent \textbf{Proof.} Set $\nu =|h|d\mu $, then the result follows from
the H\"{o}lder inequality applied to $\int fd\nu $.\hfill $\square $

We state below without proof a version of Gronwall's Lemma that plays a
crucial role in the proof of the existence and uniqueness results.

\begin{Lem}
\label{l2}(\cite{Da}, Lemma 15) Let $g:[0,T]\longrightarrow \mathbb{R}_{+}$
be a non--negative function such that $\int_{0}^{T}g(s)ds<+\infty $. Then,
there is a sequence $\{a_{n},n\in \mathbb{N}\}$ of non--negative real
numbers such that for all $p\geq 1$, $\sum_{n=1}^{+\infty }\left(
a_{n}\right) ^{1/p}<+\infty $, and having the following property: Let $%
(f_{n},n\in \mathbb{N})$ be a sequence of non--negative functions on $[0,T]$%
, $k_{1}$, $k_{2}$ be non--negative numbers such that for $0\leq t\leq T$,
\begin{equation*}
f_{n}(t)\leq k_{1}+\int_{0}^{t}(k_{2}+f_{n-1}(s))g(t-s)ds.
\end{equation*}%
If $\sup_{0\leq s\leq T}f_{0}(s)=c$, then for $n\geq 1$,
\begin{equation*}
f_{n}(t)\leq k_{1}+(k_{1}+k_{2})\sum_{i=1}^{n-1}a_{i}+(k_{2}+c)a_{n}.
\end{equation*}%
In particular, $\sup_{n\geq 0}\sup_{0\leq t\leq T}f_{n}(t)<+\infty $.\newline
If $k_{1}=k_{2}=0$, then $\sum_{n\geq 0}\left( f_{n}(t)\right) ^{1/p}$
converges uniformly on $[0,T]$.
\end{Lem}

\begin{Lem}
\label{lap1}Let $\eta $ and $\beta $ be in $(0,1)$, if $\mathbf{(H}_{\eta
}^{\alpha })$ holds then
\begin{equation*}
\int_{0}^{T}dr\int_{\mathbb{R}^{d}}\exp \left( -2r\kappa S_{\alpha }(\xi
)\right) (S_{\alpha }(\xi ))^{2\beta }\mu (d\xi )
\end{equation*}%
is finite for all $\beta $ in $]0,(1-\eta )/2[$.
\end{Lem}

\begin{Rem}
Since the spectral measure $\mu $ is non--trivial positive tempered measure,
we can ensure that there exist positive constants $c_{1}$, $c_{2}$ and $k$
such that
\begin{equation*}
c_{1}<\int_{|\xi |<k}\mu (d\xi )<c_{2}
\end{equation*}
\end{Rem}

\noindent \textbf{Proof. }Let\text{ }$\beta \in ]0,1[$, set
\begin{eqnarray*}
I(h,x) &=&\int_{0}^{T}dr\int_{\mathbb{R}^{d}}\exp \left( -2\kappa rS_{\alpha
}(\xi )\right) \left( S_{\alpha }(\xi )\right) ^{2\beta }\mu (d\xi ) \\
&=&I_{1}+I_{2}
\end{eqnarray*}%
where
\begin{equation*}
I_{1}:=\int_{0}^{T}dr\int_{\left\{ S_{\alpha }(\xi )\leq 1\right\} }\exp
\left( -2\kappa rS_{\alpha }(\xi )\right) \left( S_{\alpha }(\xi )\right)
^{2\beta }\mu (d\xi )
\end{equation*}%
and
\begin{eqnarray*}
I_{2} &:&=\int_{\left\{ S_{\alpha }(\xi )\geq 1\right\} }\left[
\int_{0}^{T}ds\exp \left( -2\kappa rS_{\alpha }(\xi )\right) \right] \left(
S_{\alpha }(\xi )\right) ^{2\beta }\mu (d\xi ) \\
&=&\int_{\left\{ S_{\alpha }(\xi )\geq 1\right\} }\left[ \frac{1-\exp
(-2T\kappa S_{\alpha }(\xi ))}{2\kappa }\right] \left( S_{\alpha }(\xi
)\right) ^{2\beta -1}\mu (d\xi ).
\end{eqnarray*}%
Clearly $I_{1}$ is finite. And the term $I_{2}$ is bounded by
\begin{eqnarray*}
\frac{1}{\kappa }\int_{\left\{ S_{\alpha }(\xi )\geq 1\right\} }\left(
S_{\alpha }(\xi )\right) ^{2\beta -1}\mu (d\xi )\leq \frac{1}{\kappa }%
\int_{\left\{ S_{\alpha }(\xi )\geq 1\right\} }\left( S_{\alpha }(\xi
)\right) ^{2\beta -1+\eta }\frac{\mu (d\xi )}{\left( 1+S_{\alpha }(\xi
)\right) ^{\eta }} && \\
\leq \frac{1}{\kappa }\sup_{\left\{ S_{\alpha }(\xi )\geq 1\right\} }\left(
S_{\alpha }(\xi )\right) ^{2\beta -1+\eta }\int_{\mathbb{R}^{d}}\frac{\mu
(d\xi )}{\left( 1+S_{\alpha }(\xi )\right) ^{\eta }}. &&
\end{eqnarray*}%
Now, choose $\beta $ in $]0,(1-\eta )/2[$ and by the hypothesis $\mathbf{(H}%
_{\eta }^{\alpha })$ the last term is finite.

\begin{Lem}
\label{lap2}If $\mathbf{(H}_{\eta }^{\alpha })$ holds then%
\begin{equation*}
\sup_{0\leq t\leq T}\int_{t}^{t+h}ds\int_{\mathbb{R}^{d}}\left\vert \mathcal{%
F}\mathbf{G}_{\alpha ,\delta }(t+h-s,x-\cdot )(\xi )\right\vert ^{2}\mu
(d\xi )\leq c(h+h^{\beta })
\end{equation*}%
for any $\beta \in ]0,1-\eta [$.
\end{Lem}

\noindent \textbf{Proof. }Let $h$ be in $[0,T]$ such that $t+h$ belongs to
be in $[0,T]$ and set
\begin{equation*}
I_{2}(h)=\int_{0}^{h}ds\int_{\mathbb{R}^{d}}\left\vert \mathcal{F}\mathbf{G}%
_{\alpha ,\delta }(s,x-\cdot )(\xi )\right\vert ^{2}\mu (d\xi ).
\end{equation*}%
By the definition of $\mathbf{G}_{\alpha ,\delta }$ we recall that $\mathcal{%
F}\mathbf{G}_{\alpha ,\delta }(s,x-\cdot )(\xi )=\psi _{\alpha ,\xi }(s)$
where $\psi _{\alpha ,\xi }(s)$ is given by (\ref{eqpsi}). Using Fubini's
theorem we can write
\begin{eqnarray*}
I_{2}(h) &=&\int_{\mathbb{R}^{d}}\left[ \int_{0}^{h}\left\vert \psi _{\alpha
,\xi }(s)\right\vert ^{2}ds\right] \mu (d\xi ) \\
&\leq &\int_{\mathbb{R}^{d}}\left[ \int_{0}^{h}\exp \left( -2\kappa
rS_{\alpha }(\xi )\right) dr\right] \mu (d\xi ) \\
&=&I_{2.1}(h)+I_{2.2}(h),
\end{eqnarray*}%
where
\begin{equation*}
I_{2.1}(h):=\int_{\left\{ S_{\alpha }(\xi )\leq 1\right\} }\left[
\int_{0}^{h}\exp \left( -2\kappa rS_{\alpha }(\xi )\right) dr\right] \mu
(d\xi )
\end{equation*}%
and
\begin{eqnarray*}
I_{2.2}(h) &:&=\int_{\left\{ S_{\alpha }(\xi )\geq 1\right\} }\left[
\int_{0}^{h}\exp \left( -2\kappa rS_{\alpha }(\xi )\right) dr\right] \mu
(d\xi ) \\
&=&\int_{\left\{ S_{\alpha }(\xi )\geq 1\right\} }\left[ \frac{1-\exp
(-2\kappa hS_{\alpha }(\xi ))}{2\kappa S_{\alpha }(\xi )}\right] \mu (d\xi ).
\end{eqnarray*}%
It is easy to see that $I_{2.1}(h)$ is bounded by $c\cdot h$. \newline
Moreover for \text{each }$\beta \in ]0,1[$ we have
\begin{equation*}
\left\vert 1-\exp \left( -2\kappa hS_{\alpha }(\xi )\right) \right\vert \leq
2^{1-\beta }\left\vert 1-\exp \left( -2\kappa hS_{\alpha }(\xi )\right)
\right\vert ^{\beta }.
\end{equation*}%
Using the mean value theorem we get
\begin{equation*}
I_{2.2}(h)\leq c_{\beta }h^{\beta }\int_{\left\{ S_{\alpha }(\xi )\geq
1\right\} }\left( S_{\alpha }(\xi )\right) ^{\beta -1}\mu (d\xi ).
\end{equation*}%
Choosing $\beta \in ]0,1-\eta [$ and use hypothesis $\mathbf{(H}%
_{\eta }^{\alpha })$ we show easily that
\begin{equation*}
\int_{\left\{ S_{\alpha }(\xi )\geq 1\right\} }\left( S_{\alpha }(\xi
)\right) ^{\beta -1}\mu (d\xi )<+\infty .
\end{equation*}

\begin{Lem}
\label{lap3}Let $z$ belongs to a compact subset $K$ of $\mathbb{R}^{d}$.
Under $\mathbf{(H}_{\eta }^{\alpha })$ we have
\begin{equation*}
\int_{\mathbb{R}^{d}}\int_{0}^{T}dr\exp \left( -2\kappa rS_{\alpha }(\xi
)\right) \left[ 1-\cos \left( 2\pi \left\langle \xi ,z\right\rangle \right) %
\right] \mu (d\xi )\leq c(\left\Vert z\right\Vert _{\infty }+\left\Vert
z\right\Vert _{\infty }^{2\beta })
\end{equation*}%
for any $0<\beta <\min \{1/2,\alpha _{0}(1-\eta )/2\}$.
\end{Lem}

\noindent \textbf{Proof.} Let us put
\begin{eqnarray*}
J(z) &=&\int_{\mathbb{R}^{d}}\left( \int_{0}^{T}dr\exp \left( -2\kappa
rS_{\alpha }(\xi )\right) ds\right) \left[ 1-\cos \left( 2\pi \left\langle
\xi ,z\right\rangle \right) \right] \mu (d\xi ) \\
&=&J_{1}(z)+J_{2}(z),
\end{eqnarray*}%
where
\begin{eqnarray*}
J_{1}(z) &:&=\int_{0}^{T}dr\int_{\left\{ S_{\alpha }(\xi )\leq 1\right\}
}\exp \left( -2\kappa rS_{\alpha }(\xi )\right) \left[ 1-\cos \left( 2\pi
\left\langle \xi ,z\right\rangle \right) \right] \mu (d\xi ) \\
&\leq &\int_{\left\{ S_{\alpha }(\xi )\leq 1\right\} }\left[ 1-\cos \left(
2\pi \left\langle \xi ,z\right\rangle \right) \right] \mu (d\xi )
\end{eqnarray*}%
and
\begin{eqnarray*}
J_{2}(z) &:&=\int_{\left\{ S_{\alpha }(\xi )>1\right\} }\frac{1-\exp \left(
-2\kappa TS_{\alpha }(\xi )\right) }{2\kappa S_{\alpha }(\xi )}\left[ 1-\cos
\left( 2\pi \left\langle \xi ,z\right\rangle \right) \right] \mu (d\xi ) \\
&\leq &\int_{\left\{ S_{\alpha }(\xi )>1\right\} }\frac{\left[ 1-\cos \left(
2\pi \left\langle \xi ,z\right\rangle \right) \right] }{2\kappa S_{\alpha
}(\xi )}\mu (d\xi ) \\
&\leq &\frac{2^{1-2\beta }}{\kappa }\int_{\left\{ S_{\alpha }(\xi
)>1\right\} }\frac{\left\vert 1-\cos \left( 2\pi \left\langle \xi
,z\right\rangle \right) \right\vert ^{2\beta }}{S_{\alpha }(\xi )}\mu (d\xi )
\end{eqnarray*}%
for any $\beta \in ]0,1/2[$.

By the mean value theorem
\begin{equation*}
J_{1}(z)\leq c\sum_{i=1}^{d}\left\vert z_{i}\right\vert \leq c\cdot d\cdot
\left\Vert z\right\Vert _{\infty }
\end{equation*}%
and
\begin{eqnarray*}
J_{2}(z) &\leq &c_{\beta ,\kappa }\int_{\left\{ S_{\alpha }(\xi )>1\right\} }%
\frac{(\sum_{i=1}^{d}\left\vert z_{i}\right\vert |\xi _{i}|)^{2\beta }}{%
S_{\alpha }(\xi )}\mu (d\xi ) \\
&\leq &c_{\beta ,\kappa ,d}\left\Vert z\right\Vert _{\infty }^{2\beta
}\int_{\left\{ S_{\alpha }(\xi )>1\right\} }\frac{(\sum_{i=1}^{d}|\xi
_{i}|)^{2\beta }}{S_{\alpha }(\xi )}\mu (d\xi ).
\end{eqnarray*}%
If $1<\alpha _{i}\leq 2$ for $i=1,\ldots ,i_{0}$, set $\alpha
_{0,1}=\min_{1\leq i\leq i_{0}}\{\alpha _{i}\}>1$. Hence
\begin{equation}
\left( \sum_{i=1}^{i_{0}}|\xi _{i}|\right) ^{2\beta }\leq c_{\beta ,d}\left(
\sum_{i=1}^{i_{0}}|\xi _{i}|^{\alpha _{0,1}}\right) ^{\frac{2\beta }{\alpha
_{0,1}}}\leq c_{\beta ,d}\left( \sum_{i=1}^{i_{0}}|\xi _{i}|^{\alpha
_{i}}\right) ^{\frac{2\beta }{\alpha _{0,1}}},  \label{a01}
\end{equation}%
where we have used the convexity property of function $x\longmapsto
x^{\gamma }$ for $\gamma >1$ and the inequality $|x|^{\alpha _{0,1}}\leq
|x|^{\alpha _{i}}$ for $\left\vert x\right\vert >1$.\newline
If $0<\alpha _{i}<1$ for $i=i_{0},\ldots ,d$, set $\alpha
_{0,2}=\min_{i_{0}\leq i\leq d}\{\alpha _{i}\}$. Use the fact that $%
|x|^{1/\alpha _{i}}\leq |x|^{1/\alpha _{0,2}}$ for $\left\vert x\right\vert
>1$ and the super--additivity of the function $x\longmapsto x^{\gamma }$ for
$\gamma >1$ to obtain
\begin{equation}
\left( \sum_{i=i_{0}}^{d}|\xi _{i}|\right) ^{2\beta }\leq c_{\beta ,d}\left(
\sum_{i=i_{0}}^{d}|\xi _{i}|^{\alpha _{i}\frac{1}{\alpha _{0,2}}}\right)
^{2\beta }\leq c_{\beta ,d}\left( \sum_{i=i_{0}}^{d}|\xi _{i}|^{\alpha
_{i}}\right) ^{\frac{2\beta }{\alpha _{0,2}}}.  \label{a02}
\end{equation}%
Recall that $\alpha _{0}=\min_{1\leq i\leq d}\{\alpha _{i}\}$, then $\alpha
_{0}=\min \{\alpha _{0,1},\alpha _{0,2}\}$. Combining (\ref{a01}) and (\ref%
{a02}) with concavity of the mapping $x\longmapsto x^{\gamma }$ $(\gamma <1)$
for $x>1$, we end up with
\begin{equation*}
\left( \sum_{i=0}^{d}|\xi _{i}|\right) ^{2\beta }\leq c_{\beta ,d}\left(
\sum_{i=0}^{d}|\xi _{i}|^{\alpha _{i}}\right) ^{\frac{2\beta }{\alpha _{0}}},
\end{equation*}%
for $2\beta <\alpha _{0}$. Consequently
\begin{equation*}
\int_{\left\{ S_{\alpha }(\xi )>1\right\} }\frac{(\sum_{i=1}^{d}|\xi
_{i}|)^{2\beta }}{S_{\alpha }(\xi )}\mu (d\xi )\leq c_{\beta ,\kappa
,d}\int_{\left\{ S_{\alpha }(\xi )>1\right\} }\left( S_{\alpha }(\xi
)\right) ^{\frac{2\beta }{\alpha _{0}}-1}\mu (d\xi ),
\end{equation*}%
which is finite according to $\mathbf{(H}_{\eta }^{\alpha })$ for all $%
0<\beta <\alpha _{0}(1-\eta )/2.$\hfill $\square $

\end{document}